\documentclass{amsart}
\usepackage{graphicx}
\usepackage{amsthm,amsfonts,amssymb}
\usepackage{amsmath}
\usepackage{hyperref}
\usepackage[latin1]{inputenc}
\usepackage{latexsym}

\newtheorem{theo}{Theorem}[section]
\newtheorem{conj}[theo]{Conjecture}
\newtheorem{coro}[theo]{Corollary}
\newtheorem{lemm}[theo]{Lemma}
\newtheorem{prop}[theo]{Proposition}

\theoremstyle{definition}
\newtheorem{rema}[theo]{Remark}

\newcommand \vsp{\vspace{0.2cm}}

\newcommand \tend{\longrightarrow}
\newcommand \be{\begin{eqnarray*}}
\newcommand \ee{\end{eqnarray*}}
\newcommand \ben{\begin{eqnarray}}
\newcommand \een{\end{eqnarray}}

\newcommand{\bigO}[1]{\mathcal O\pa{#1}}
\newcommand{\bigtheta}[1]{\Theta\pa{#1}}

\def \build#1#2#3{\mathrel{\mathop{\kern 0pt#1}\limits_{#2}^{#3}}}

\def \esp#1#2{\mathbb E_{#1}\left[#2\right]}

\def\rompar(#1){\textup(#1\textup)}    

\newcommand{\ac}[1]{\left\{#1\right\}}
\newcommand{\pa}[1]{\left(#1\right)}
\newcommand{\cro}[1]{\left[#1\right]}
\newcommand{\abs}[1]{\left|#1\right|}
\newcommand{\pr}[1]{\mathbb P\left(#1\right)}
\newcommand{\var}[1]{\textnormal{Var}\left(#1\right)}
\newcommand{\cov}[1]{\textnormal{Cov}\left(#1\right)}

\newcommand \bel{\begin{eqnarray}\label}
\newcommand \eel{\end{eqnarray}}

\title[Marcus--Lushnikov processes and union--Find algorithms]{Merging costs for the additive Marcus--Lushnikov process, and  Union-Find  algorithms.}

\begin{document}

\author{Philippe Chassaing}
\author{R{\'e}gine Marchand}
\address{Institut Elie Cartan Nancy (math{\'e}matiques)\\
Universit{\'e} Henri Poincar{\'e} Nancy 1\\
Campus Scientifique, BP 239 \\
54506 Vandoeuvre-l{\`e}s-Nancy  Cedex France}

\email{chassain@iecn.u-nancy.fr, marchand@iecn.u-nancy.fr}


\begin{abstract}
Starting with a monodisperse configuration with $n$ size--1 particles,  an additive
 Marcus--Lushnikov process evolves until it reaches its final state (a  unique particle
with mass $n$). At each of the $n-1$ steps of its evolution, a merging  cost is incurred,
 that depends on the sizes of the two particles involved, and on an  independent random factor.
This paper studies the asymptotic behaviour of the  cumulated costs up  to the
$k$th  clustering, under various regimes for $(n,k)$, with applications  to the study
of Union--Find algorithms.
\end{abstract}

\maketitle


\section{Introduction, models and results}


Fundamental to computer science is the manipulation of
 \textit{dynamic} sets: sets that can grow, shrink or
otherwise change over time.
 Some  algorithms, e.g. Kruskal or Prim algorithms
 for the search of the minimum spanning tree of a graph,
 involve grouping $n$ distincts elements into a collection of
disjoint  sets, and implementing
two operations, U{\small NION}, that unites two sets, and F{\small IND}   that finds
which set a given element belongs to (see \cite[Part III]{CLR} for  more). For the analysis of the cost of such
 operations, Yao \cite{YAO}  suggested two models,  the spanning tree  model
 and the random graph model. Both are instances of a general  model
 of coalescence of particles, that we describe now.

\subsection{Marcus--Lushnikov processes}
\label{coal}
The study  of coalescence of particles (sets, clusters) with
different sizes has a long story, and  has applications
in many  scientific disciplines besides computer science, such as
 physical chemistry, but also astronomy,
bubble swarms, and mathematical genetics  (cf. the survey \cite{Al}).
In a basic model, clusters with different masses
move through space, and when two clusters
 (say, with masses $x$ and $y$)
are sufficiently close, there is some chance
 that they merge into a single cluster
with mass $x+y$, with a probability
 quantified, in some sense, by a \textit{rate kernel}
$K$, depending on the masses, the positions and the velocities of the  two clusters.
  However, such a model, including the spatial
distribution of clusters and their velocity, is still too
 complicated for analysis, so a rather natural first approximation
was suggested independently by  Marcus  \cite{marcus}
  and Lushnikov \cite{lush1,lush2},
 by considering kernels depending only
  on the masses of the  clusters.

A \textit{Marcus--Lushnikov process}
 \cite{Al} with rate $K$
is a continuous-time Markov process whose state space is the set
 of partitions of  $n$ or, equivalently,
 the set of  measures on the set $\mathbb{N}$ of positive integers
 \[
\mu
=
\sum_k\ \frac{n(k,t)}n\ \delta_k,
\]
in which $n(k,t)$ is an integer, and
\[
\sum_kkn(k,t)=n,
\]
 so that $\int x \mu(dx)=1$.
The $k$'s stand for the sizes of clusters and $n(k,t)$ is the number
 of clusters with size $k$ at time $t$. The size--$k$ clusters provide
 a fraction $\frac{k\,n(k,t)}n$ of the total size $n$.
A Marcus--Lushnikov process
 evolves by instantaneous jumps according to the rule
\[\textrm{each pair }(x,y)\textrm{ of
 clusters  merge at rate }K(x,y)/n.\]
In other words, the system of clusters
 jumps from the state $\mu$ to the state
$\mu+\frac1n\pa{\delta_{x+y}- \delta_x-\delta_y}$ at rate $K(x,y)/n$,
meaning  that, if at time $t$ the state of the
 system is $(x_i)_{i\ge 1}$, the next pair $(I,J)$
of clusters that  merge and the time $t+T$ when they merge
 are jointly distributed as follows:
assume we are given a set of independent
 random variables $(T_{i,j})_{1\le i<j}$
with exponential distribution described by
\[\mathbb{P}\left(T_{i,j}>t\right)=\exp\left(-K(x_i,x_j)t/n\right),\]
and set
\[
\inf_{1\le i<j}T_{i,j}
=
T_{I,J}=T.
\]
It follows, as usual for continuous time
 Markov chains,  that $T_{I,J}$  and $(I,J)$
are independent, that $T_{I,J}$ has an exponential law with
 parameter $\sum_{i,j}K(x_i,x_j)$, and that
\begin{equation}
\label{pmerge}
\mathbb{P}\pa{(I,J)=(i,j)}=\frac{K(x_i,x_j)}
{\sum_{k,\ell}K(x_k,x_{\ell})}.
\end{equation}

We shall see later that the \textit{additive Marcus--Lushnikov process}
(with kernel $K(x,y)=x+y$) is embedded in  the spanning tree
 model of Yao.
 The relation between the random graph model and
the \textit{multiplicative
 Marcus--Lushnikov process} (with kernel $K(x,y)=xy$)
was noted by Knuth and Sch{\"o}nhage
\cite{KS}  and Stepanov \cite{Stepanov}.
In both cases, the
clusters are connected components of a graph,
and the merging  of two clusters is due to the addition of an
 edge between  elements of these clusters. Also,
we assume that the initial state
 consists in $n$ clusters with size 1; this state is often
 called the \textit{monodisperse configuration}.
This corresponds to a totally disconnected
graph with $n$ vertices  and no edges.
Thus there are eventually  $n-1$ jumps (steps, mergings \dots)
between the initial state $\delta_1$ and the final state $\frac1n\  \delta_n$
of the Marcus--Lushnikov process. In this paper, we focus on the additive
case.

\subsection{Analysis of merging costs}
\label{section:mergingcosts}
At the $k$-th jump (addition of the $k$-th edge)
of the Marcus--Lushnikov process, two subsets  with respective
 sizes $\pa{S_{k,n}, s_{k,n}}$,  $S_{k,n}\ge s_{k,n}$,
 are merged,
at a cost $c_{k,n}$ that may depend
on the sizes   $\pa{S_{k,n}, s_{k,n}}$.
For instance, in some implementations, a label
is maintained for each element,
signaling the set it belongs to,  and when merging
two sets,
one has to change the labels of the elements of one of the 2 sets.
Yao, Knuth and Sch{\"o}nhage
 studied two   algorithms:
\begin{itemize}
\item \textit{Quick-Find}, that updates the labels of
one of the two sets, selected arbitrarily,
leading to cumulated costs
\[C^{QF}_{n,m}=\sum_{k=1}^{m}A_{k,n},\]
in which
$A_{k,n}=S_{k,n}$ with probability $1/2$ and  $A_{k,n}=s_{k,n}$ with  probability $1/2$,
\item and
 \textit{Quick-Find-Weighted}, that updates the smaller set at a cost
$c_{k,n}=s_{k,n}$,
leading to  cumulated costs
\[C^{QFW}_{n,m}=\sum_{k=1}^{m}s_{k,n}.\]
\end{itemize}

In other contexts where coalescence of
 two sets occurs,
 costs of interest are $L_{k,n}$, the size of one of the two sets
chosen randomly with a probability that is proportional to its size,  i.e.
$L_{k,n}=S_{k,n}$ with probability $S_{k,n}/(S_{k,n}+s_{k,n})$ and
$L_{k,n}=s_{k,n}$ with probability $s_{k,n}/(S_{k,n}+s_{k,n})$, or
\begin{equation}
\label{defprey}
R_{k,n}
=
S_{k,n}+s_{k,n}-L_{k,n},
\end{equation}
or again
\[D_{k,n}= \left\lfloor U_k L_{k,n}\right\rfloor.\]
In the next Sections, some   interpretations are given for these last  costs.
Here,  $\pa{U_k}_{1\le k\le n-1} $  denotes a sequence
of independent random variables, uniform on $[0,1]$.

In \cite{KS},  using recurrence
relations, Knuth and Sch{\"o}nage give the following
 equivalents for the total merging costs:
\begin{equation}
\label{res_conv}
\esp{}{C^{QF}_{n,n-1}}
=
\sqrt{\frac{\pi}{8}}\ n^{3/2}+O(n \log n),
\hspace{0,5cm}
\esp{}{C^{QFW}_{n,n-1}}
=
 \frac{1}{\pi}n \log n + O(n),
\end{equation}
in the  case of the  additive Marcus--Lushnikov process
($\log$ denotes  the natural logarithm).
In this paper,  we study concentration or limit laws for total costs
$C_{n,n-1}$ as well as for partial costs
$C_{n,\lceil\alpha n\rceil}$. For the partial costs,
we obtain the following results:

\begin{theo}
\label{coutpart_biased}
For any $\eta\in (0,1)$, and any positive $\varepsilon$,
$$
\lim_n\pr{\sup_{\alpha \in [0,1-\eta]}
\left|
 \frac{C^{QF}_{n,\lceil\alpha n\rceil}}{n} - \varphi^{QF}(\alpha)
\right|\ge\varepsilon}
=
0,
$$
respectively
$$
\lim_n\pr{\sup_{\alpha \in [0,1-\eta]}
\left|
 \frac{C^{QFW}_{n,\lceil\alpha n\rceil}}{n} - \varphi^{QFW}(\alpha)
\right|\ge\varepsilon}
=
0,
$$
in which
\begin{eqnarray*}
\varphi^{QF}(\alpha)
& = &
\frac12 \left( \frac{1}{1-\alpha} + \log \left(\frac{1}{1 -\alpha}  \right)
\right), \\
\varphi^{QFW}(\alpha)
& = &
\int_0^{ \log \left(
     \frac{1}{1 -\alpha} \right)} \sum_{k \in \mathbb{N}} \sum_{l \in
  \mathbb{N}} (k\vee l)\  q(k,t)q(l,t)\ dt, \\
q(k,t)
& = & 
\frac{\left[k(1-e^{-t})\right]^{k-1} e^{-t}}{k!}\  \exp(-k(1-e^{-t})).
\end{eqnarray*}
\end{theo}
This Theorem is actually a corollary
of Theorem \ref{coutpart_general}.
Theorem \ref{coutpart_general} is stated and proven at Section
\ref{smolu}:
it gives the expression, in terms of the solution  $q(k,t)$ of the
Smoluchowski equation, of the limit
 function $\varphi^c(\alpha)$ for the partial costs:
\[
C_{n,\lceil\alpha n\rceil}
=
\sum_{k=1}^{\lceil\alpha n\rceil}\hat c\pa{S_{k,n},s_{k,n},U_{k,n}}
\]
once $C_{n,\lceil\alpha n\rceil}$ is normalized by $\frac1n$.
For Theorem \ref{coutpart_general} to cover a wide class of costs (starting with Quick Find),
the general expression  $\hat c\pa{S_{k,n},s_{k,n},U_{k,n}}$ for the instantaneous cost
of the $n$--th jump
has to involve
an extra--randomization parameter, $U_{k,n}$,  uniform on $[0,1]$.
Theorem \ref{coutpart_general} holds true under the mild condition
of polynomial growth, as a function of $S_{k,n}$ and $s_{k,n}$,
of the instantaneous conditional cost
\[
c(S_{k,n},s_{k,n})
=
\esp{}{
\left.\hat c\pa{S_{k,n},s_{k,n},U_{k,n}}
\right\vert\pa{S_{k,n},s_{k,n}}
}.
\]
For instance, the instantaneous conditional cost for Quick Find is
\[
\esp{}{
\left.A_{k,n}
\right\vert\pa{S_{k,n},s_{k,n}}
}
=
\frac{S_{k,n}+s_{k,n}}2.
\]

For QFW and QF,
 the total costs are respectively  $\bigtheta{n\log n}$ or  $\bigtheta{n^{3/2}}$,
 while the partial costs are $\bigtheta{n}$:
this is consistent with
\[
\lim_{1}\varphi^c(\alpha)
=
+\infty,
\]
and also, of course,
$
\esp{}{C^{QFW}_{n,n-1}}
=
o\pa{\esp{}{C^{QF}_{n,n-1}}}
$
is consistent with
$\varphi^{QFW}=o\pa{\varphi^{QF}}$.
Note that, compared with \cite{KS}, Theorem \ref{coutpart_biased}
adds some kind of concentration result for partial costs. We turn
now to a more precise study of the total costs.

\subsubsection*{Detailed analysis of the total cost for QFB and
QFW}
Let us define
\[
C^{QFB}_{n,m}
=
\sum_{k=1}^{m}R_{k,n}.
\]
An interpretation of $R_{k,n}$ in terms of the spanning tree model is given
 in the next Sections
(QFB stands for Quick-Find-Biased).
 We have
\begin{theo}
\label{olebotheo_right}
$$
\frac{C^{QFB}_{n,n-1}}{n \log n}
\ \build{\longrightarrow }{}{\mathcal L_2}\
\frac 12.
$$
\end{theo}
{From} \eqref{defprey},
$R_{k,n}=S_{k,n}$ with probability $s_{k,n}/(S_{k,n}+s_{k,n})$ and
$R_{k,n}=s_{k,n}$ with probability $S_{k,n}/(S_{k,n}+s_{k,n})$.
As a consequence
 $R_{k,n}$ is more likely equal to
the smaller block $s_{k,n}$ than to  $S_{k,n}$, so we expect   similar  behaviours for
$C^{QFB}_{n,n-1}$ and  $C^{QFW}_{n,n-1}$. Moreover we
expect a smaller variance for  $C^{QFW}_{n,n-1}$ than for
$C^{QFB}_{n,n-1}$,
but we could not produce a proof.
However, at the light of Theorem
\ref{olebotheo_right}, we conjecture that
\begin{conj}
\label{olebotheo_min}
$$
\frac{C^{QFW}_{n,n-1}}{n \log n}
\ \build{\longrightarrow }{}{\mathcal L_2}\
\frac 1\pi.
$$
\end{conj}

\subsubsection*{Detailed analysis of the total cost for Quick-Find}
Let $(e(t))_{0 \leq t \leq 1}$ denote
the normalized Brownian excursion. For $C^{QF}_{n,n-1}$,
we have the following result:

\begin{theo}\label{olebotheo_random}
$n^{-3/2}\,C^{QF}_{n,n-1}$ converges in law to $\int_0^1
e(t)dt.$
\end{theo}
Actually, a more precise result is available: for $\beta\ge 0$, let
\begin{eqnarray*}
W_n(\beta)
&=&
n^{-3/2}C^{QF}_{n,\lfloor n- \beta\sqrt n\rfloor}
\\
&=&
n^{-3/2}\sum_{k=1}^{\lfloor n- \beta\sqrt n\rfloor}A_{k,n},
\\
h_{\beta}(t)
&=&
e(t)-\beta t-\inf_{0\le s\le  t}\left(e(s)-\beta s\right),\\
W(\beta)
& = &
\int_{0}^{1} h_{\beta}(t) dt.
\end{eqnarray*}
Then
\begin{theo}
\label{process}
$\pa{W_n(\beta)}_{\beta\ge 0}$ converges in law to
$\pa{W(\beta)}_{\beta\ge 0}.$
\end{theo}
\noindent Theorem \ref{olebotheo_random}
is the convergence of $W_n(0)$.
For a detailed study of the family
$\pa{W(\beta)}_{\beta\ge 0}$, see \cite{Ja}.
Since  $\lim_{+\infty}W(\beta)=0$, Theorem \ref{process}
 yields that:
\begin{coro}
\label{phasetransition}
Assume that $\sqrt n=o(h_n)$ and $h_n\le n$. Then
$$
n^{-3/2}C^{QF}_{n,\lfloor n-h_n\rfloor}
\build{\longrightarrow }{}{P}
0.
$$
\end{coro}

\begin{rema}
\label{rem:comparaison}
As  opposed to  Quick--Find,
the partial sums for Quick--Find--Biased
satisfy
$$
\lim_{n}\pa{n \log n}^{-1}\esp{}{C^{QFB}_{n,\lfloor n-h_n\rfloor}}
=
\lim_{n}\pa{n \log n}^{-1}\esp{}{C^{QFB}_{n, n-1}},
$$
for $h_n=o(n)$, and the same property holds for Quick--Find--Weighted.  These
 quite different behaviours for the partial and total costs of QF and  QFW
can be explained, partly, by the existence of several
different regimes of convergence
of the additive Marcus--Lushnikov process.
\end{rema}

\subsection{Regimes of the additive Marcus--Lushnikov process.}
Denote by $B^n_{k,1}$ the size of the
largest cluster after the $k$--th jump:
interpretations based on fragmentation of trees \cite {ADD,PAV}
or on analysis of hashing algorithms \cite {CL}
show that the additive Marcus--Lushnikov process
has  three different regimes:

\begin{itemize}
\item
the \textit{sparse regime}: if $\sqrt n =o(n-k)$, then $B^n_{k,1}/n  \rightarrow 0$
in probability ;
\item
the  \textit{transition regime}:
 when $n-k=O(\sqrt n)$,  several clusters of size $O(n)$
coexist, and, once renormalized, clusters' sizes converge to the
 widths of excursions of
Brownian-like stochastic processes  ;
\item
the \textit{almost full  regime}:
if $n-k=o(\sqrt n)$, $B^n_{k,1}/n \rightarrow 1$ in probability,
and a unique giant cluster of size $n-o(n)$
 coexists with smallest clusters with total size $o(n)$.
\end{itemize}

Thus,  the dramatic increase
of  $B^n_{k,1}$
(and, as a consequence, of $A_{k,n}$) during
the transition regime explains
the huge contribution of the transition regime
to  the sum  $C^{QF}_{n,n-1}$, as quantified by
Theorem \ref{process}
and by Corollary \ref{phasetransition},
and this in spite of the fact that the transition regime
 involves a relatively small number of  terms of  $C^{QF}_{n,n-1}$.
Rather than $B^n_{k,1}$,
 the sizes of small clusters have an actual impact on
$C^{QFW}_{n,n-1}$ or $C^{QFB}_{n,n-1}$,
since, in most of the jumps,
$s_{k,n}$ is way smaller than $S_{k,n}$ ;
thus  the quite different behaviour of QF and QFB reveals that,
 in some sense,
the sizes of small clusters have a   moderate increase during the  transition regime,
the sparse
regime providing the largest contribution
to $C^{QFW}_{n,n-1}$ or $C^{QFB}_{n,n-1}$.
Also, the apparition of the Brownian excursion
 area in Theorems \ref{olebotheo_random} and
\ref{process} is typical of a phenomenon
linked with the transition regime, 
where the asymptotics of the parking scheme
can be described in terms of the standard additive coalescent
 \cite{ADD, BERT, CL}.

The asymptotic behaviour of the partial costs
$C_{n,\lfloor \alpha n\rfloor}$
is determined by
the behaviour of the additive Marcus--Lushnikov
process during the sparse regime:
once suitably normalized, the additive Marcus--Lushnikov process  converges to the
(deterministic) solution of Smoluchowski equations
(cf. \cite{Fou_Gie,norris} or Theorem \ref{th:smolu--MLush}),
 explaining the deterministic nature of the limits
$\varphi^{QF}(\alpha)$ and $\varphi_{QFW}(\alpha)$ in Theorem  \ref{coutpart_biased}.

The paper is organized as follows:
in Section \ref{embeddings}, we describe the embedding of the additive  Marcus--Lushnikov process
in two  combinatorial coalescence
models, the random spanning tree and the parking scheme.
Through the first embedding, we can rephrase the analysis of Union-Find algorithms in terms of the
additive 
Marcus--Lushnikov process. Convergence of Marcus--Lushnikov processes to
solutions of Smoluchowski equations  is used in
Section \ref{smolu} to prove Theorem \ref{coutpart_biased}.
In Section \ref{sectionproofQFW}, we use some combinatorial properties  of
the parking scheme to bound the mean and the variance of
Quick-Find-Biased and  prove Theorem
\ref{olebotheo_right}.
In
Sections \ref{QFfullsection}
 and \ref{almostfullsection},
 we prove  Theorems
\ref{olebotheo_random} and \ref{process} about the total cost of  Quick-Find,
with the help of the analysis of phase transitions for the parking, as  given in
\cite{CL}.

\section{Two embeddings of the additive Marcus--Lushnikov process}
\label{embeddings}

Marcus--Lushnikov processes  are of no use to
 Knuth, Sch{\"o}nhage or Yao,
and  their analysis of average costs
 of UNION-FIND algorithms rely quite naturally on probabilistic models
defined in terms of  random spanning trees,
 or in terms of random graphs.
Following \cite{Pi},
 the next subsection  recalls  how the additive
 Marcus--Lushnikov process  $X^{(n)}
 =\pa{X^{(n)}_t}_{t\ge 0}$ is embedded
in the spanning tree model.
As a consequence, the analysis of partial costs
for the additive
 Marcus--Lushnikov process, given in Section
 \ref{smolu},
  turns out to be a development of
  Knuth, Sch{\"o}nhage or Yao analysis.
The proofs  of Sections \ref{sectionproofQFW}--\ref{almostfullsection}
 rely on the embedding of the additive
 Marcus--Lushnikov process
in  the parking model, a model often used to analyze
 linear probing in hashing tables \cite{CL,FPV}.
This last embedding is described in a second subsection.

We start with a description of the additive Marcus--Lushnikov process
that helps to understand its connections to the spanning tree model and
 to the parking scheme:
at  step $k$ pick a  first cluster  $P$ with a probability
$\frac{|P|}{n}$ among the $n-k+1$ clusters, and let us call it the  ``predator"
(being  a size--biased pick it is likely larger than the average  cluster) ;
then pick the ``prey" $p$ uniformly among the $n-k$
remaining clusters, and let $P$ eat $p$, producing a unique cluster  with size
$|P|+|p|$. It is not hard to see that  this defines
the additive Marcus--Lushnikov process, and  that
 $L_{k,n}$  (resp. $R_{k,n}$)
can be seen as the size of the predator (resp.
of the prey).  If, alternatively, both clusters are size--biased picks
 (resp. if both are uniform picks), we
obtain the \textit{multiplicative Marcus--Lushnikov process}
 (resp. the \textit{Marcus--Lushnikov process with constant kernel}, also
 called Kingman's process).

\subsection{The spanning tree model.}

 Let $\mathcal T_n$ be the set of unrooted  labeled trees
with $n$ vertices. As noted by Cayley, $\mathcal T_n$
 has  $n^{n-2}$ elements.
Given a labeled tree $T\in \mathcal T_n$,
 consider a  labelling (or ordering) of its $n-1$ edges.
Let $T_k$  be the subgraph of $T$ whose $k$  edges have labels
 not larger than $k$: $T_k$ is a forest with
$n-k$ connected components. The connected components  (trees)
of the forest play the role
of the  dynamic sets we mentioned earlier. We have:
\begin{itemize}
\item
$T_0$ is the graph with no edges. It  has  $n$
  size-1  components, that
 we call \textit{monomeres}, following
 chemists' terminology. Also, $T_{n-1}=T$.
\item
$T_k$ is obtained from $T_{k-1}$ by addition
of the edge labelled $k$ in $T$.
\end{itemize}

Following \cite{KS}, let us call the sequence $(T_k)_{0 \leq k \leq  n-1}$ a
\textit{spanning tree of} $T$. Now, there are $(n-1)!$ orderings
of the $n-1$ edges of this
tree, and thus the set $ST_n$ of
spanning trees has $n^{n-2}\times (n-1)!$  elements.
A \textit{random spanning tree} is a random uniform element of $ST_n$.

Let $Y_k$ be the partition of the  number $n$ induced by the
 connected components of $T_k$.  In  \cite{Pi}, Pitman proves that
 conditionally
given  $(Y_i)_{0 \leq i \leq k}$, the addition of the $k+1$-th edge will
 merge two subtrees with respective sizes $x$ and
$y$    with a probability
\[\frac{x+y}{n(n-k-1)}.\]
The same expression is obtained  specializing relation (\ref{pmerge})
 to the case  $K(x,y)=a(x+y)$, when $X^{(n)}_t$ has exactly $k$  clusters.
Thus  $Y^{(n)}=(Y_i)_{0 \leq i \leq n-1}$ and
$X^{(n)}=\pa{X^{(n)}_t}_{t\ge 0}$ have the
 same law,  up to a time change: the jumps of $Y^{(n)}$
take place at times 1, 2, \dots ,  $n$, while  the jumps of $X^{(n)}$
occur at random times
\footnote
{
However an exact
identity between the two processes
is easily  obtained through a standard randomization artifice:
attach independent exponential random times $t_e$ with
mean 1 to each edge $e$ of  a random uniform labeled tree
$T\in \mathcal T_n$, and let  the edge $e$ appear at time $t_e$.
Let $T_t$ be the subgraph of $T$ with edges $e$ such that $t_e\le t$,
and let $Y^{(n)}_t$ be the partition of $n$ induced by the
connected components of $T_t$. Then
 $Y^{(n)}=\left( Y^{(n)}_t \right)_{t\ge 0}$ is
a Marcus--Lushnikov process with kernel $K(x,y)=(x+y)/n$.
}
 (actually the time elapsed
between the $k$-th and  $k+1-$th jumps  of $X^{(n)}$ is
random exponentially distributed with mean
$\frac1{an(n-k-1)}$).
As the merging costs do not depend on the precise times
of  jumps, but only on the sizes of clusters that merge,
this difference does not matter:  the total and partial costs  have
the same law  in  the additive Marcus--Lushnikov process and
 in the spanning tree model. Thus the Yao--Knuth--Sch{\"o}nhage
problem fits in the more general frame
of merging costs for   Marcus--Lushnikov processes.

In this context, $R_{k,n}$ and  $L_{k,n}$  have the following
 interpretation:  let  any fixed vertex  be the root, once and for all,
so that each
edge has a bottom vertex (the vertex that is closer to the root)
and a top vertex.
Erasing the $k$--th edge splits a subtree  of
$T_k$  in two connected components (clusters),
the  ordered sizes of our clusters being $s_{k,n}\le S_{k,n}$,
with the notations of Section \ref{section:mergingcosts}.
It turns out that  the size of the cluster at the bottom of the
$k$--th  edge is a size--biased pick among
$\ac{s_{k,n}, S_{k,n}}$. Thus   $L_{k,n}$  (resp. $R_{k,n}$)
can be seen as the size of the cluster at the bottom (resp.
at the top) of the $k$--th edge, just before the $k$--th
jump.

\subsection{The parking model.}
\label{subsection:parking}
Consider a parking lot of $n$ places on a roundabout, on
which a set $\mathcal C=\{1,2,\dots,n-1\}$ of $n-1$
cars eventually park. Each car $c$ has a clock that rings
at a time $T_c$, and when the clock rings,
the car $c$ tries to park on a random  place $t(c)$.
If the first try $t(c)$ is on an empty place, the car parks there;
otherwise, the car tries the next places clockwise,
and parks on the first empty place it finds.
The first tries $\pa{t(c)}_{c\in \mathcal C}$
 are assumed independent
 and uniform on the $n$ places, numbered from 1 to $n$, and
times $\pa{T_c}_{c\in \mathcal C}$ are assumed
to be independent exponentially  distributed, with mean 1.

In this model, the clusters are the blocks of places
 already occupied, with  the following conventions:
\begin{itemize}
\item there are as many  blocks as there are empty places,
\item a block contains an empty place and the set of consecutive
 occupied   places before (going clockwise) this empty place,
\item the size of the block is the total number of places in it,
 including the empty place,
\item if an empty place follows another empty place,
it is considered as a size--1 block of its own.
\end{itemize}

\begin{figure}[tbp]
\includegraphics[width=4.6cm]{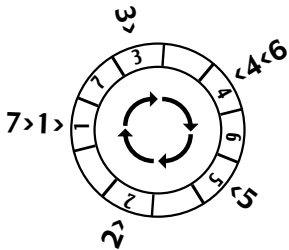}
\caption{A sample of tries $t(c)$ and the resulting 3 clusters.}
Here
$n=10=4+4+2$.
\end{figure}

This way, the initial configuration, with $n$ empty places,
has $n$   size--1 blocks (i.e. is monodisperse), and each time a car  parks,
two blocks merge, with conservation of the mass, as the empty place
 that disappears and the car that
replaces it both count for one mass--unit.
The final configuration,
once the $n-1$ cars are parked, has a unique cluster with size $n$,
and a unique empty place, with number $V$ uniformly distributed on
$\{1,2,\dots,n\}$.

It turns out that the sizes of blocks
form an additive Marcus--Lushnikov process, with kernel  $K(x,y)=(x+y)/n$:
given that the parking scheme
 with  $n$ places, $k$ cars already parked
and  $\ell=n-k$ empty places,
has two blocks with sizes $x\ge y$, the probability that
 these two blocks merge at the next arrival is
\begin{equation}
\label{pt}
\frac{x+y}{n(n-k-1)}.
\end{equation}

Actually, as follows from equiprobability for the $n^{k}$ possible
 configurations,
 the number $N_{x,y}$ of empty places after block
$x$ (clockwise) but before block $y$
 is random uniform on $1,2,\dots, \ell-1$.
If $N_{x,y}\notin \{1,\ell-1\}$,  there is no way  the two clusters
can merge  at the next arrival.
Given that  $N_{x,y}=1$ (resp. $\ell-1$)
the conditional probability
 that the two blocks merge at the next arrival is the probability
that the next time a clock ring, the first try of the corresponding car
 will be on one of the $x$ (resp. $y$) places of the largest (resp.  smallest)
cluster:
\[
\frac{x}{n},
\hspace{0.5cm}\mbox{resp.}\hspace{0.2cm}
\frac{y}{n},
\]
leading to \eqref{pt}.
Another consequence is that the size  of the block before
 (clockwise) the place filled by the
$k$--th arrival  is a random size--biased choice among
$\ac{s_{k,n}, S_{k,n}}$:
 $L_{k,n}$ and $R_{k,n}$ can be seen  as
the sizes of  blocks before (clockwise)  and after the place filled by the
$k$--th arrival, and $D_{k,n}$ as the displacement of the car between
its first try and its final place.

\vspace{0.2cm}
{From} the parking interpretation, we deduce now some explicit  computations
for the law of the weighted blocks $L_{k,n}$ and
$R_{k,n}$, that give some light on the asymptotic behaviour of
 $S_{k,n}$ and
$s_{k,n}$.
Consider the conditional probability $p_{m,k}^{(j,n)}$
that, in an additive Marcus--Lushnikov process
with size $n$, the $j$-th predator
has size $k$, before  the $j$-th meal, given that its size after the  $j$-th meal is $m$.
{From} now on, we assume  the Marcus--Lushnikov process
 to be embedded  in a parking scheme. In particular, we retain
 the interpretation of $L_{k,n}$ and $R_{k,n}$ as
the sizes of  blocks before and after the place filled by the
$k$--th arrival, so that $p_{m,k}^{(j,n)}$
is the probability that, in a parking  scheme with $n$ places,
  the  block before the place  filled (resp. the block created)
 by the $j$-th arrival
has size $k$ (resp. $m$).
It turns out, for combinatorial reasons,
 that  $p_{m,k}^{(j,n)}$
does not depend
on $j$ or $n$. Thus  we have, for instance,
$$
p_{m,k}^{(j,n)}
=
p_{m,k}^{(m-1,m)}
=
\mathbb{P}(L_{m-1,m}=k)
=
\mathbb{P}(R_{m-1,m}=m-k),
$$
and  we shall  drop the exponent, for seek of brevity.
{From} the asymptotic behaviour of
 $p_{m,k}$, we expect
some  intuition about the respective values of
$L_{k,n}$ and $R_{k,n}$.

\begin{lemm}
\label{lemm:comb-basic}
$$p_{m,k}=\frac{1}{m^{m-2}} {{m-2}\choose{k-1}}
k^{k-1} (m-k)^{m-k-2}.$$
\end{lemm}

\begin{proof}
Recall that  the size of a  cluster
is defined as the number of cars in the block \textit{plus one}.
There are ${{m-2}\choose{k-1}}$ possible choices  for the
$k-1$ cars in the block after  $V$ (clockwise) ,
and $k^{k-2}$ possible parking schemes
for these cars ; also, there are $(m-k)^{m-k-2}$
possible parking schemes
for  the $m-k-1$  cars in the block before  $V$,
and finally, $k$ possible first tries for the last car if  $V$
is to be the last empty place.
\end{proof}

\vsp
Lemma \ref{lemm:comb-basic} and Stirling's formula yield at once that

\begin{coro}
\label{rem_loideR}
\begin{equation} \label{convenloi}
\forall k \geq 1, \lim_{m \rightarrow \infty}
p_{m,m-k}=\frac{k^{k-1}e^{-k}}{k!}.
\end{equation}
\end{coro}

The limit distribution is the so--called
Borel  distribution, tightly related
to explicit solutions of Smoluchowski equations \cite{Al}, and
to the tree function or Lambert's function \cite{KNU2}.
Thus, in distribution,
$R_{m-1,m}=\bigO{1}$ in some sense.
However, note  that the
Borel  distribution has infinite mean,
 in coherence with the fact that
$\esp{}{R_{m-1,m}}=\Theta\pa{\sqrt m}$.
We shall retain that, provided
$L_{k,n}
+
R_{k,n}
$
is large,
$R_{k,n}$ or $s_{k,n}$  are negligible, compared with
$L_{k,n}$. As a consequence,
 $S_{k,n}$ or  $L_{k,n}$ should have quite similar behaviours.
 This is a first tentative explanation of the drastic difference between
 QF and QFW, revealed by Knuth \& Sch{\"o}nhage' results.

 \begin{rema}
 \label{tagged}
The convergence of the Marcus--Lushnikov process to the solution
of the Smoluchowski equation, derived by analytic arguments in
\cite{norris}, is quite natural for the additive case at the light of the
following computations.
The probability $p(\alpha n)$ that, after the $\alpha n$--th arrival,
 the first car to be parked belongs to a size--$k$
cluster, is
\[
\frac
{{\alpha n\choose k-2}k^{k-2}(n-k)^{\alpha n-k}(n-\alpha n-1)n}
{n^{\alpha n}}
\sim
(1-\alpha)\alpha^{k-2}\ \frac{k^{k-2}}{(k-2)!}\ e^{-\alpha k}.
\]
As  the size--$k$ clusters provide a fraction $\frac{kn(k,t)}{n}$ of  the total size, they also provide
a fraction $\frac{(k-1)n(k,t)}{n(t)}$ of the total number $n(t)$ of cars
arrived at time $t$, so the probability $\hat p(t)$ that, at time $t$,
 the first car to be parked belongs to a size--$k$
cluster is precisely  $\frac{(k-1)n(k,t)}{n(t)}$. We shall see later that the
$\alpha n$--th arrival takes place at a  time $t_\alpha\sim -\log(1-\alpha)$, so that 
$\hat p( -\log(1-\alpha))\sim p(\alpha n)$,
or, equivalently:
\[
\frac{(k-1)}\alpha \ \frac{n(k, -\log(1-\alpha))}{n}
\sim
(1-\alpha)\alpha^{k-2}\ \frac{k^{k-2}}{(k-2)!}\ e^{-\alpha k}.
\]
The right hand side turns out to be the expression 
of   $\frac{k-1}{\alpha}\ q(k,-\log(1-\alpha))$.
 \end{rema}

\section{Analysis of partial costs after $\lceil\alpha n\rceil$  coalescences}
\label{smolu}
In this Section we state and prove Theorem \ref{coutpart_general}, and
  Theorem \ref{coutpart_biased} follows as a direct consequence. As
opposed to the next Sections,
the proofs make no use of richer combinatorial structures
in which the additive Marcus--Lushnikov process is embedded,
and they could very likely be generalized
to a suitable class of kernels $K$.
We assume that the cost incurred at the $k$th step is
\[
\tilde{c}_{k,n}
=
\hat{c}(s_{k,n},S_{k,n},U_{k,n})\geq 0
\]
in which $(U_{k,n})_{k \in \mathbb{N},
n\in \mathbb{N}}$ denote a sequence of
independent identically distributed random variables uniform on $[0,1]$:
 this covers the case of QFW, in which
the cost $A_{k,n}$ can be written
\[
A_{k,n}
=
s_{k,n}\mathbf{1}_{U_{k,n}\le 0.5}+S_{k,n}\mathbf{1}_{U_{k,n}>0.5}.
\]
The size of the prey $L_{k,n}$ can be written
\[
L_{k,n}
=
s_{k,n}\mathbf{1}_{U_{k,n}\le \frac{S_{k,n}}{s_{k,n}+S_{k,n}}}
+S_{k,n}\mathbf{1}_{U_{k,n}> \frac{S_{k,n}}{s_{k,n}+S_{k,n}}},
\]
the size of the predator and the displacement
have similar descriptions.
We suppose that there exist $A>0$ and  $p,q \in
\mathbb{N}$ such that:
$$\forall x \in \mathbb N, \; \forall y \in \mathbb N, \;  h(x,y)=\int_0^1
\hat{c}^2(x,y,u)du \leq Ax^py^q.$$
We set, for $1\le m\le n-1$,
$$
C_{n,m}
=
\sum_{k=1}^{m}\tilde{c}_{k,n}.
$$
Then the asymptotic behaviour of
$
C_{n,\lceil\alpha n\rceil}
$
can be described in terms  of
the instantaneous conditional cost
\begin{eqnarray*}
c(x,y)&=&\int_0^1 \hat{c}(x,y,u) du
\\
&=&\esp{}{\left.\hat c\pa{S_{k,n},s_{k,n},U_{k,n}}
\right\vert\pa{S_{k,n},s_{k,n}}=(x,y)},
\end{eqnarray*}
and of the solution of  the Smoluchowski equation
with additive kernel (see Subsection \ref{smolu2} below):
\[
q(k,t)
=
\frac{\left[k(1-e^{-t})\right]^{k-1} e^{-t}}{k!}\  \exp(-k(1-e^{-t})).
\]
We have
\begin{theo}
\label{coutpart_general}
For any $\eta>0$,
$$\sup_{\alpha \in [0,1-\eta]}
\left|
 \frac{C_{n,\lceil\alpha n\rceil}}{n} - \varphi^c(\alpha)
\right|
\ \build{\longrightarrow }{}{P}\ 0,$$
in which $\varphi^c$ is an increasing function
from $[0,1)$ to ${\mathbb{R}}^+$ defined by
\[
\varphi^{c}(\alpha)
=
\int_0^{ \log \left(
     \frac{1}{1 -\alpha} \right)} \sum_{k \in \mathbb{N}} \sum_{l \in
  \mathbb{N}} c(k,l)\  q(k,t)q(l,t)\ dt.
\]
\end{theo}
Thus, $\varphi^c$ corresponds to a renormalized partial cost  until time $\log \left(
     \frac{1}{1 -\alpha} \right)$ in the
infinite particle system governed by Smoluchowski equation.
In the table below, we give  the explicit values of $\varphi^c$ for
some examples: 

\vspace{0.4cm}
\begin{center}
\noindent
\begin{tabular}{|l|c|c|}
\hline
Cost & $c(x,y)$ & $\varphi^c(\alpha)$ \\
\hline
\hline
Quick-Find $A_{k,n}$
&
$ \displaystyle \frac{x+y}{2}$
&
$ \displaystyle \frac12 \left( \frac{1}{1-\alpha} + \log  \left(\frac{1}{1 -\alpha} \right) \right)$
\\
\hline
Prey size $L_{k,n}$
&
$\displaystyle  \frac{2xy}{x+y}$
&
$\log \left(\frac{1}{1 -\alpha} \right)$
 \\
\hline
Predator size $R_{k,n}$
&
$ \displaystyle \frac{x2+y2}{x+y} $
&
$ \displaystyle  \frac{1}{1 -\alpha}$
\\
\hline
Displacement $D_{k,n}$
&
$ \displaystyle \frac{x2+y2}{2(x+y)} $
&
$ \displaystyle \frac{1}{2(1 -\alpha)}$
\\
\hline
\end{tabular}
\end{center}

\vspace{0.2cm}
\noindent
For  Quick-Find-Weighted, $c(x,y)$ has the simple form $\min(x,y)$,
but we could not produce an expression more explicit than
\[
\varphi^{QFW}(\alpha)
=
\int_0^{ \log \left(
     \frac{1}{1 -\alpha} \right)} \sum_{k \in \mathbb{N}} \sum_{l \in
  \mathbb{N}} (k\vee l)\  q(k,t)q(l,t)\ dt.
\]

Note that  a similar expression appears in the analysis of Union-Find
algorithms under the random graph model (kernel
$K(x,y)=xy$): 
  Bollob{\'a}s \&  Simon \cite{BOLL}
proved that
the average cost of QFW  is $cn+O(n/\log n)$, in which:
\[
c
=
\log
2-1+\sum_{k\ge1}
\pa{
\frac1{k}-\frac{k^k}{k!}
\sum_{\ell=1}^{k-1}\frac{\ell^{\ell-1}}{\ell!}\frac{k+\ell -2!}{(k+\ell)^{k+\ell-1}}
}.
\]

\subsection{The additive Smoluchowski equation.}
\label{smolu2}
The  proof of Theorem \ref{coutpart_general} relies on the
convergence of the additive Marcus--Lushnikov process to
the solution of  the Smoluchowski equation
with additive kernel. Let   $\mathcal M_1^+(\mathbb{N})$ denote the set of
positive measures on $\mathbb{N}$ with total mass less or equal to 1.
A (deterministic) solution $\mu$ of the additive Smoluchowski
equation is a family
$\mu=(\mu_t)_{t \geq 0}$ of measures in $\mathcal M_1^+(\mathbb{N})$
\[
\mu_t
=
\sum_{k \in \mathbb{N}}q(k,t)\delta_k,
\]
that satisfy:
$$(S) \left\{
\begin{array}{ll}
i) & \forall k \in \mathbb{N}, \; q(k,0)=\delta_1(k), \\
ii) & \forall k \in \mathbb{N}, \; \forall t \geq 0, \\
&\\
& \; \; \frac{dq(k,t)}{dt}=\frac{1}{2} \sum_{j=1}^{k-1}
kq(j,t)q(k-j,t)-q(k,t)\sum_{j=1}^\infty (j+k)q(j,t).
\end{array}
\right.
$$
The coefficient $q(k,t)$ can be seen as the concentration of particles
of size $k$ at time $t$ in a given volume unit, for an
infinite system of particles. The first term on the right hand side
of the Smoluchowski equation  $(S)$ corresponds to the creation of
a particle with size $k$ due to coalescence between smaller particles,
of size $j$ and $k-j$, at a rate $j+(k-j)=k$,
and the second term to the destruction of a particle with size $k$,
through coalescence with another particle of size $j$, at a rate $k+j$.

In the additive case, there exists a unique solution to $(S)$, given by:
\begin{equation}
\label{solution+}
\forall k \in \mathbb{N}, \; \forall t \geq0, \; q(k,t)=\frac{1}{k}
\frac{\left[k(1-e^{-t})\right]^{k-1}}{(k-1)!} e^{-t-k(1-e^{-t})}
\end{equation}
(see
Aldous \cite{Al}). All the moments of this solution can be
explicitly computed, and for instance:
$$\forall t\geq 0, \; <\mu_t,x>=1, \;
<\mu_t,1>=e^{-t}, \; <\mu_t,x2>=e^{2t}.$$
The first equality says that the mass is preserved during coalescences,
the second one says that the concentration (number of particles per unit
volume) decreases exponentially, and the third one gives the exponential
increase of the mean size of a tagged (size biased) particle.

\subsection{The infinitesimal generator of the additive
Marcus--Lushnikov process.}
An alternative definition of the additive Marcus--Lushnikov process,
through its infinitesimal
generator, is more suitable for our computations.
An additive
 Marcus--Lushnikov
process $(\mu^n_t)_{t \geq 0}$ is a continuous time c{\`a}dl{\`a}g
Markov process
with values in $\mathcal M_1^+(\mathbb{N})$,
satisfying the set $(ML_{n})$ of conditions below:
\begin{itemize}
\item[i.] $\mu^n_0=\delta_1,$
\item[ii.] $ \forall t \geq 0, \mu^n_t \in \{ \frac{1}{n} \sum_{i=1}^k
\delta_{x_i}, \; k \in \mathbb{N}, \; \forall i \; x_i\in \mathbb{N}, \;
\sum_{i=1}^k x_i =n\},$
\item[iii.]
its generator $L$ is given by:
$$
\left\{\begin{array}{ll}
& \forall \psi:M_{1}^+(\mathbb{N}) \rightarrow \mathbb{R}
 \mbox{ measurable}, \;
 \forall \mu=\frac{1}{n} \sum_{i=1}^k \delta_{x_i}, \\
&\\
& \; \; L\psi(\mu)=\sum_{i \neq j}
\left(\psi(\mu+\frac{1}{n}(\delta_{x_i+x_j}-\delta_{x_i}  -\delta_{x_j}))-\psi(\mu)\right)\left( \frac{x_i+x_j}{2n}\right).
\end{array}
\right.
$$

\end{itemize}
In the last term, for
symetry reasons, the additive kernel appears with  a factor 1/2.

It is well known
 that, for every $n$,  $(ML_n)$ has a unique solution
$(\mu^n_t)_{t \geq 0}$ (which is a collection of
 random measures in $\mathcal M_1^+(\mathbb{N})$), satisfying moreover
to the mass conservation property: $$
\forall t\geq 0, \; <\mu^n_t,x>=1 \; \; a.s.$$

\subsection{Convergence of the solution of $(ML_n)$ to the solution of  $(S)$.}
We recall here some definitions and theorems  of convergence for the  additive
Marcus--Lushnikov process.

\vspace{0.2cm}
\noindent
1. On $\mathcal M_1^+(\mathbb{N})$,  the vague  convergence of measures
is defined as follows:
$$(\mu_n)_{n \in \mathbb N} \ \build{\longrightarrow }{}{v}\ \mu
\; \Leftrightarrow \; \forall
\psi \in C_c(\mathbb N, \mathbb{R}), \; <\mu_n, \psi> \rightarrow
<\mu,\psi>,$$
in which $C_c(\mathbb N, \mathbb{R})$ denotes the space of functions  from
$\mathbb N$ to $\mathbb{R}$ with compact support. We assume that
$\mathcal M_1^+(\mathbb{N})$ is endowed with the vague topology (which
is metrizable).
Denote by
$\mathbb{D}([0,T], \mathcal M_1^+(\mathbb{N}))$ the set of
c{\`a}dl{\`a}g functions from $[0,T]$ to $\mathcal M_1^+(\mathbb{N})$,
endowed with the Skorokhod topology \cite{Eth_Kur}.

Denote by $(\mu_t^n)_{t \geq 0}$ the
solution of $(ML_n)$ and by $(\mu_t)_{t \geq 0}$ the solution of $(S)$.
Our analysis makes use of the following convergence theorem (it is a refinement, due to \cite{Fou_Gie},
of a well known result of \cite{norris}), and of some
direct consequences listed below:
\begin{theo}
\label{th:smolu--MLush}
For every $T>0$,
$$(\mu_t^n)_{t \in [0,T]}  \ \build{\longrightarrow }{}{dist}\  (\mu_t)_{t \in
  [0,T]}.$$
\end{theo}
\noindent Here we mean convergence in distribution.

\vspace{0.2cm}
\noindent
2. As $(\mu_t)_{t \geq 0}$ is deterministic, the  convergence in
distribution implies the convergence in probability, that is, if $d$  denotes a
metric yielding the Skorokhod topology on $\mathbb{D}([0,T], \mathcal
M_{1}^+(\mathbb{N}))$, we have:
$$\forall T>0, \; \forall \varepsilon>0, \hspace{1cm}\mathbb{P}\pa{
d\cro{(\mu_t^n)_{t \in [0,T]},(\mu_t)_{t \in
  [0,T]}}\geq \varepsilon} \longrightarrow 0.$$

\vspace{0.2cm}
\noindent
3. Since the limit $t \mapsto \mu_t$ is continuous,
convergence for the Skorokhod topology
entails  uniform convergence on
every $[0,T]$: for any metric  $d_v$   yielding the vague topology on
$M_{\leq 1}^+(\mathbb{N})$, we have
$$\forall T>0, \; \forall \varepsilon>0,  \hspace{1cm} \mathbb{P} \left( \sup_{t  \in
    [0,T]} d_v[\mu^n_t,\mu_t] \geq \varepsilon \right) \longrightarrow  0.$$

\vspace{0.2cm}
\noindent
4. Finally,  we have
\begin{prop}
\label{polyn}
For any
function
$\varphi$ from
$\mathbb{N}$ to $\mathbb{R}$ satisfying, for some $A>0$ and $p \in
\mathbb{N}$, $|\varphi(k)| \leq Ak^p$,
$$\forall T>0, \; \forall \varepsilon>0,  \hspace{1cm} \mathbb{P} \left( \sup_{t  \in
    [0,T]} |<\mu^n_t,\varphi>-<\mu_t,\varphi>|\geq \varepsilon  \right)
  \longrightarrow 0.$$
\end{prop}
\noindent When  $\varphi$ is a function from
$\mathbb{N}$ to $\mathbb{R}$ with compact support,
Proposition \ref{polyn} follows directly from point 3,
but for the  class of functions with
 polynomial growth,
we  need some bounds
 on the moments $<\mu_t,x^p>$ and
$\esp{}{<\mu^n_t,x^p>}$:

\begin{lemm}
\label{momentsMLS}
For every $p \geq 2$, there exist positive constants $A_p$
  and $B_p$ such that for every $t \geq 0$:
\begin{eqnarray}
\label{momentsMLS2}
\esp{}{<\mu^n_t,x^p>}
&\leq&
e^{B_pt},\\
\label{momentsMLS1}
<\mu_t,x^p>
&\leq&
A_p e^{2(p-1)t}.
\end{eqnarray}
\end{lemm}

\begin{proof}
We derive relation \eqref{momentsMLS2}
 using the special  form of the infinitesimal
generator of a Marcus--Lushnikov process
(cf. $(ML_n)$). To this aim, some additional notations are handy:
for a function $\psi$
from $\mathbb{N}^2$ in $\mathbb{R}^+$ and a measure
$\mu=\frac{1}{n} \sum_{i=1}^k \delta_{x_i}$,
let us define
$$
<\mu \stackrel{\Delta_n}{\otimes} \mu, \psi>
=
<\mu \otimes \mu, \psi>\ -\ \frac1n\int\psi(x,x)\mu(dx).
$$
When $\mu=\frac{1}{n} \sum_{i=1}^k \delta_{x_i}$,
then
$$
<\mu \stackrel{\Delta_n}{\otimes} \mu, \psi>
=
\frac{1}{n2} \sum_{i\neq j} \psi(x_i,x_j).
$$
We have
$$\esp{}{<\mu^n_t,x^p>}=1+\int_0^t \esp{}{<\mu^n_s
\stackrel{\Delta_n}{\otimes} \mu^n_s, \left( (x+y)^p-x^p-y^p \right)
\left( \frac{x+y}{2} \right)>}ds.$$
Since $\left( (x+y)^p-x^p-y^p \right)
\left( \frac{x+y}{2} \right) \leq (2^{p-1}-1)(x^py+y^px)$, for all $x$  and $y$ in
$[0,+\infty)$,
\begin{eqnarray*}
\esp{}{<\mu^n_t,x^p>}
& \leq &
1+(2^{p-1}-1) \int_0^t \esp{}{<\mu^n_s
\stackrel{\Delta_n}{\otimes} \mu^n_s, x^py+y^px>}ds
\\
& \leq &
1+(2^{p-1}-1)\int_0^t \esp{}{<\mu^n_s{\otimes} \mu^n_s,
x^py+y^px>}ds
\\
& \leq & 1+(2^{p}-2) \int_0^t \esp{}{<\mu^n_s,x^p>}ds,
\end{eqnarray*}
the last relation making use of the mass
conservation property.
Now \eqref{momentsMLS2} follows from Gronwall's Lemma.
Similar technics lead to inequality \eqref{momentsMLS1}, the complete proof
can be found in  \cite{Dea_Tan}.
\end{proof}

\begin{proof}[Proof of Proposition \ref{polyn}.]
We consider
\begin{eqnarray*}
\alpha_{K,n}
&=&
\mathbb{P} \left( \sup_{t \in [0,T]} |<\mu^n_t-\mu_t,\varphi\  \mathbf{1}_{[0,K)}>|
\geq \varepsilon/3  \right),
\\
\beta_{K}
&=&
\sup_{t \in [0,T]} |<\mu_t,\varphi\ \mathbf{1}_{[K,+\infty)}>|,
\\
\gamma_{K,n}
&=&
\mathbb{P} \left( \sup_{t \in [0,T]} |<\mu^n_t,\varphi\  \mathbf{1}_{[K,+\infty)}>|
\geq \varepsilon/3  \right).
\end{eqnarray*}
First,
\begin{eqnarray*}
\esp{}{\sup_{t \in [0,T]} |<\mu^n_t,\varphi\ \mathbf{1}_{[K,+\infty)}>|}
&\le&
A\,\esp{}{\sup_{t \in [0,T]} <\mu^n_t,x^p\ \mathbf{1}_{[K,+\infty)}>}
\\
&\le&
A\,K^{-p}\,\esp{}{\sup_{t \in [0,T]} <\mu^n_t,x^{2p}>}
\\
&\le&
A\,K^{-p}\,\esp{}{<\mu^n_T,x^{2p}>},
\end{eqnarray*}
the last inequality due to the fact that
$t\,\rightarrow\,<\mu^n_t,x^{2p}>$  is increasing,
as a consequence of $a^p+b^p\le(a+b)^p$.
Thus \eqref{momentsMLS2} and Markov inequality lead to a uniform bound
\[
\gamma_{K,n}
\le
3\ A\,K^{-p}\,e^{B_pT}\varepsilon^{-1}.
\]
Also,
\[
\beta_K
\le
A_{2p} \,K^{-p}\,e^{2(2p-1)T}.
\]
As a consequence, $K$ can  be tuned to make
$\sup_n \gamma_{K,n}$
arbitrary small, and  simultaneously $\beta_K$  smaller than  $\varepsilon /3$.
 Once $K$  chosen, we use $\lim_n \alpha_{K,n}=0$ to conclude.
\end{proof}

\vspace{0.2cm}
\noindent
5. By a similar proof,  for every function $\psi$ from
$\mathbb{N}^2$ to $\mathbb{R}$ such that $ |\psi(k,l)| \leq  Ak^pl^q,$  we have 
\begin{equation}
\label{convf2}
\lim_n \mathbb{P} \left( \sup_{t \in
    [0,T]} |<\mu^n_t\otimes \mu^n_t,\psi>-<\mu_t \otimes  \mu_t,\psi>|\geq \varepsilon  \right)
=
0,
\end{equation}
for any $T$ and $\varepsilon$ positive.

\subsection{Merging costs as functionals of $(ML_n)$.}
In this subsection, we prove  Theorem \ref{coutpart_general}.  
Let $(U^n_s)_{s \geq 0}$ denote a family of  independent and
identically distributed random variables, uniform on $[0,1]$ and independent
of $(\mu_{t}^n)_{t\ge 0}$. When a coalescence
occurs at time $s$  ($\mu_{s-}^n\neq
\mu_s^n$),  we assume that  a nonnegative cost $\tilde{c}(\mu_{s-}^n,
\mu_s^n,U^n_s)$ is incurred, with 
$$
\tilde{c}\left(
\frac{1}{n}\sum_{i=1}^k \delta_{x_i}, \frac{1}{n}\sum_{i=1}^k  \delta_{x_i} +
\frac{1}{n} \left( \delta_{x_i+x_j}-\delta_{x_i}-\delta_{x_j}\right),u 
\right)
=
\hat{c}(x_i,x_j,u)
$$
if $k \in \{2,\dots,n\}, \;  (x_i)_{1 \leq i
  \leq k} \in \mathbb{N}^k, $ and $ u \in [0,1]$,
  and with
  $
  \tilde{c}\left(
\mu,\nu,u 
\right)
$
null otherwise. Furthermore, we assume that there exist $A>0$ and  $p,q \in
\mathbb{N}$ such that:
$$
h(x,y)
=
\int_0^1
\hat{c}^2(x,y,u)du
\leq
 Ax^py^q,
\hspace{1cm}
\forall x \in \mathbb N, \; \forall y \in \mathbb N.
$$
Then the partial cost up to time $t$ is
$$C^n_t  =  \sum_{0<s\leq t} \tilde{c}(\mu_{s-}^n, \mu_s^n,U^n_s).$$
Recall that $c(x,y)=\int_0^1 \hat{c}(x,y,u)du$. According to
\cite[Ch. IV, Lemma (21.13)]{Rog_Wil}, we have
\begin{eqnarray*}
\frac{C^n_t}{n}
& = & \int_0^t <\mu^n_s
\stackrel{\Delta_n}{\otimes}  \mu^n_s, c(x,y)\frac{x+y}2>ds +M^n_t
\\
& = &
\int_0^t <\mu^n_s \otimes \mu^n_s, c(x,y)\frac{x+y}2>ds -
\frac{1}{n} \int_0^t <\mu^n_s,xc(x,x)>ds +M^n_t,
\end{eqnarray*}
in which $M^n_t$ is a martingale such that
$$<M^n>_t=\frac{1}{n} \int_0^t <\mu^n_s
\stackrel{\Delta_n}{\otimes}  \mu^n_s, h(x,y)\frac{x+y}{2}>ds.
$$
Set
\begin{eqnarray*}
C_t
&=&
\int_0^t <\mu_s \otimes \mu_s, c(x,y)\frac{x+y}2>ds
\\
&=&
\int_0^t \int \int c(x,y)d \mu_s(x)d\mu_s(y)ds.
\end{eqnarray*}
As a consequence of  the convergence of the solution $(\mu_t^n)_{t \geq 0}$  of $(ML_n)$ to the
solution $(\mu_t)_{t \geq 0}$ of $(S)$,
we get:
\begin{theo}
\label{coutMLS}
For every cost $\hat{c}$ such that there exist $A>0$ and  $p,q \in
\mathbb{N}$ with \\
$\forall x \in \mathbb N, \; \forall y \in \mathbb N, \; h(x,y)=\int_0^1
\hat{c}^2(x,y,u)du \leq Ax^py^q$, we have, for each positive $T$ and  $\varepsilon$,
\[
\lim_n\mathbb{P} \left( \sup_{t \in [0,T]} \left| \frac{C^n_t}{n}-C_t  \right|
  \geq \varepsilon \right)
=
0.
\]
\end{theo}

\begin{proof}\  First we bound the martingale and  the diagonal term.
By Doob's inequality, we obtain
\begin{eqnarray}
\label{martin}
\esp{}{\sup_{t \in [0,T]} |M^n_t|}^2
&\leq&
4\,\esp{}{<M^n>_T},
\end{eqnarray}
but Lemma \ref{momentsMLS}  yields that
\begin{eqnarray*}
\esp{}{<M^n>_t}
& \leq &
\frac{A}{2n} \int_0^t \esp{}{<\mu^n_s \otimes \mu^n_s,x^py^q(x+y)>}ds
\\
& \leq &
\frac{A}{2n} \int_0^t
\left( e^{(B_{p+1}+B_q)s}+e^{(B_{p}+B_{q+1})s} \right)ds,
\end{eqnarray*}
that vanishes as $n$ grows to infinity.
For the diagonal term
$$D^n_t=\frac{1}{n} \int_0^t <\mu^n_s,xc(x,x)>ds,$$
observe that
\begin{eqnarray*}
D^n_t\leq \frac{\sqrt A}{n} \int_0^t <\mu^n_s,x^{p+q+1}>ds,
\end{eqnarray*}
and that $t\rightarrow D^n_t$ is increasing.
Thus it is enough to
control the terminal value:
\begin{eqnarray}
\nonumber
\esp{}{ \sup_{t \in [0,T]} D^n_t}
&\le&
\frac{\sqrt A}{n} \int_0^T \esp{}{<\mu^n_s,x^{p+q+1}>}ds
\\
\label{diag}
&\leq&
\frac{\sqrt A\,T\,e^{B_{p+q+1}T}}{n},
\end{eqnarray}
that vanishes as $n$ grows to infinity.
Then, with the help of  (\ref{convf2}),  we bound the integral terms:  for any positive
$ T$ and $ \varepsilon$, we have
$$
\lim_n\mathbb{P} \left( \sup_{t \in [0,T]} \left|
\int_0^t <\mu^n_s \otimes \mu^n_s-\mu_s
\otimes \mu_s, c(x,y)\frac{x+y}2>ds \right|\geq \varepsilon \right)
=
0.
$$
Finally, as usual,
\begin{eqnarray*}
\lefteqn{\mathbb{P} \left( \sup_{t \in [0,T]} \left|  \frac{C^n_t}{n}-C_t \right|
\geq \varepsilon \right)}
\\
& \leq &
\mathbb{P} \left( \sup_{t \in [0,T]} \left|
\int_0^t <\mu^n_s \otimes \mu^n_s-\mu_s\otimes \mu_s,  c(x,y)\frac{x+y}2>ds
\right|\geq \varepsilon/3 \right)
\\
&  &\hspace{2,6cm}+ \mathbb{P} \left( \sup_{t \in [0,T]}  \left|M^n_t\right| \geq
  \varepsilon/3 \right)
 + \mathbb{P} \left( \sup_{t \in [0,T]} D^n_t\geq
  \varepsilon/3 \right),
\end{eqnarray*}
and the three terms on the right hand side vanish,
 the first one by  step 2, the second (resp.  third) term,
by   \eqref{martin} (resp. \eqref{diag}) and by Markov inequality.
\end{proof}

\begin{proof}[Proof of Theorem \ref{coutpart_general}.]
For  analysis of algorithms or combinatorics, the fact that
Marcus--Lushnikov processes are \textit{continuous--time} processes
looks like an artefact:
 this artefact will prove  useful if we can convert Theorem  \ref{coutMLS},
a result about the cumulated cost at
a deterministic \textit{time}, into a result  about the cumulated cost  after
a deterministic \textit{number of jumps}.
Thus we have to establish a close connection between
the cumulated cost $C^n_{t}$ up to time $t$,
defined at the previous section, and the cumulated costs $C_{n,m}$ or  $C_{n,\lceil\alpha n\rceil}$
 involved in  Theorem \ref{coutpart_general}.
For $\alpha \in [0,1)$, set:
\[
T^n_\alpha
=
\inf\ac{t\geq 0, \;<\mu^n_t,1>\ \leq \ 1-\alpha-\frac1n};
\]
$T^n_\alpha$ is the time when the $\lceil\alpha n\rceil$-th
coalescence occurs, when  the total number of clusters becomes
smaller than $(1-\alpha)n-1$. Thus
\begin{equation}
\label{equn:mu}
<\mu^n_{T^n_\alpha},1>
\ \simeq\
1-\alpha,
\end{equation}
and
\begin{equation}
\label{equn:cost}
C^n_{T^n_\alpha}
=
C_{n,\lceil\alpha n\rceil}.
\end{equation}
 As a consequence of Proposition \ref{polyn},
for any positive $T$ and $\varepsilon$, we have
\[
\lim_n\mathbb{P} \left( \sup_{t \in [0,T]} |<\mu^n_t,1>-<\mu_t,1>|\geq  \varepsilon  \right)
=
0.
\]
Since $<\mu_t,1>=e^{-t}$, relation \eqref{equn:mu}  leads
to $e^{-T^n_\alpha}\sim 1-\alpha$, and
the following Lemma is not unexpected:
\begin{lemm}
\label{TT}
For any positive  $\varepsilon$ and $\eta$,
\[
\lim_n
\mathbb{P} \left( \sup_{\alpha \in [0,1-\eta]} \left| T^n_\alpha+
\log \left( 1-\alpha\right) \right| \geq\varepsilon   \right)
=
0.
\]
\end{lemm}

\begin{proof}
\ Assume that for some $\alpha \in [0,1-\eta]$, we have:
\[
T^n_\alpha+\log \left( 1-\alpha\right)
\geq
\varepsilon,
\]
or
\[
T^n_\alpha+\log \left( 1-\alpha\right)
\leq
-\varepsilon.
\]
The first inequality insures that for any time
$t_0<-\log \left( 1-\alpha\right)+\varepsilon\le \varepsilon-\log \eta$,
$<\mu^n_{t_0},1>$ is larger than $1-\alpha$, and if for instance we  choose
$t_0>-\log \left( 1-\alpha\right)+\varepsilon/2$, we obtain
\[
\abs{<\mu^n_{t_0},1>-<\mu_{t_0},1>}
>
\eta(1-e^{-\varepsilon/2}).
\]
The second inequality
insures that at time
\[
t_1
=
-\log \left( 1-\alpha\right)-\varepsilon\ge 0,
\]
we have
$
<\mu^n_{t_1},1>\
\le
\ 1-\alpha,
$
and as a consequence
\[
\abs{<\mu^n_{t_1},1>-<\mu_{t_1},1>}
>
\eta(1-e^{-\varepsilon/2}).
\]
Then  we use Proposition \ref{polyn}, with $T=\varepsilon-\log \eta$.
\end{proof}

\vspace{0.2cm}
Finally, we   combine relation \eqref{equn:cost}, Theorem \ref{coutMLS}  and
 Lemma \ref{TT} to deduce the proof of Theorem \ref{coutpart_general}.
Recall that
\begin{eqnarray*}
\varphi(\alpha)
&=&
C_{ \log {\left( \frac{1}{1-\alpha} \right)} }.
\end{eqnarray*}
Given any positive numbers $\beta$, $\varepsilon$ and $\eta$, we can  write:
\begin{eqnarray*}
\lefteqn{\mathbb{P} \left( \sup_{\alpha \le1-\eta} \left|
    \frac{C^n_{T^n_\alpha}}{n} -C_{ \log {\left( \frac{1}{1-\alpha}
  \right)} } \right| \geq \varepsilon \right)}
\\
& \leq & \mathbb{P} \left( \sup_{\alpha \le1-\eta} \left| T^n_\alpha+
    \log (1-\alpha) \right| \geq \beta \right)\\
&& + \mathbb{P} \left( \sup_{\alpha \le1-\eta} \left|
    \frac{C^n_{T^n_\alpha}}n - C_{T^n_\alpha}\right| \geq  \varepsilon/2, \;
  \sup_{\alpha \le1-\eta} \left| T^n_\alpha +
    \log (1-\alpha)\right| \leq \beta \right) \\
&& + \mathbb{P} \left( \sup_{\alpha \le1-\eta} \left| C_{T^n_\alpha} -
    C_{ \log {\left( \frac{1}{1-\alpha} \right)} } \right| \geq  \varepsilon/2
, \;\sup_{\alpha \le1-\eta} \left| T^n_\alpha +
    \log (1-\alpha) \right| \leq \beta \right)  \\
 & \leq & \mathbb{P} \left( \sup_{\alpha \le1-\eta} \left| T^n_\alpha +
    \log ( 1-\alpha) \right| \geq \beta \right)
 +  \mathbb{P} \left( \sup_{t \le \beta-\log \eta}
 \left| \frac{C^n_{t}}{n} - C_{t}\right| \geq \varepsilon/2 \right) \\
&&\hspace{4,4cm}+  \mathbf{1}_{\left\{ \sup\ac{ \left. \left|  C_{t}-C_{s}  \right|\ \right\vert \
s,t \in \left[ 0, \beta -\log \eta  \right],\; |t-s| \leq \beta }
 \geq \varepsilon /2\right\}}.
\end{eqnarray*}
For  $\beta$ small enough
 the third term of the last sum vanishes,
by the  uniform continuity of $t \mapsto C_t$.
Theorem \ref{coutMLS} and
 Lemma \ref{TT} take care of the two other terms.
\end{proof}

\section{Analysis of the total cost of Quick-Find-Biased}
\label{sectionproofQFW}

\subsection{Average case analysis} 
 In this subsection, as a first step for the proof of Theorem
\ref{olebotheo_right},
 we prove the convergence of the
first moment of $R_n/n \log n $, using the parking representation. 
In  the next subsection,  a bound
 for the variance of  $R_n/n \log n $ completes
the proof of Theorem
\ref{olebotheo_right}. We have:
\label{section:Meancost}
\begin{lemm}
\label{conv_en_moy_droit}
\[
\lim_n\ \frac{\esp{}{C^{QFB}_{n,n-1}}}{n \log n }
=
\frac{1}{2}.
\]\end{lemm}

The next Lemma is of constant use in the rest of the paper:

\begin{lemm}
\label{equationRL}
For any $k\in\{1, \dots, n-1\}$,
$\displaystyle
\esp{}{R_{k,n}|L_{k,n} }=\frac{n-L_{k,n}}{n-k}.
$
\end{lemm}
\begin{proof}
As in Section \ref{subsection:parking}, we assume  the  Marcus--Lushnikov process
 to be embedded  in a parking scheme.
 Let us number the blocks clockwise from 0 to $n-k$,
starting with the block  before  the place filled by the
$k$--th arrival, and let $\beta_i$ denote the size of the $i$--th block
(so that $(\beta_0,\beta_1)=(L_{k,n},R_{k,n})$).
It is easy to see that
 among the $n^k$ parking configurations, there are
\begin{equation}
\label{costadditivity}
{k-1\choose b_0-1,b_1-1,\dots,b_{n-k}-1}n  b_0\prod_{i=0}^{n-k}b_i^{b_i-2}
\end{equation}
configurations  such that
$\pa{\beta_i}_{0\le i\le n-k}=\pa{b_i}_{0\le i\le n-k}$.
As a consequence,  the family $\pa{\beta_i}_{1\le i\le n-k}$ is
 exchangeable,
while $\beta_0$, being a size--biased pick among the $n-k+1$ blocks,
tends to be larger. With the additional fact that
\[
\sum_{i=0}^{n-k}\beta_i
=
n,
\]
this leads to
\[
\esp{}{\beta_i\left|\beta_0\right.}=\frac{n-\beta_0}{n-k},
\]
for any $i\ge 1$, and specially for $\beta_1=R_{k,n}$.
\end{proof}

\begin{proof}[Proof of Lemma \ref{conv_en_moy_droit}.]
We find different bounds for $\esp{}{R_{k,n} }$ according to the three different
regimes of the additive Marcus--Lushnikov process.
 For $\varepsilon$  positive but smaller than $1/2$,  set
$\varphi(n)=n-n^{\frac{1}{2}+\varepsilon}$ and
$\psi(n)=n-n^{\frac{1}{2}-\varepsilon}$. Also, let
  $B^n_{k,1}\ge B^n_{k,2}\ge  \dots$ denote the sequence of sizes of
  blocks (clusters) after the $k$--th
arrival (jump), in decreasing order:
\subsubsection{The sparse regime.}
For $k\leq\varphi(n)$,  the largest cluster is small,
and, as a consequence,
$$\esp{}{R_{k,n} }=\frac{n-\esp{}{L_{k,n} }}{n-k}\simeq \frac{n}{n-k},$$
or, more precisely,
\begin{lemm} \label{coucoumevoila}
$\displaystyle
\lim_n
\ \ \sup\ac{\left.
1-\frac{n-k}{n}\ \esp{}{R_{k,n} }\right|\ 1 \leq k \leq \varphi(n)}= 0.
$
\end{lemm}
\begin{proof}
By  Lemma \ref{equationRL},
\begin{equation}
\label{aaa}
1-\frac{n-k}{n}\ \esp{}{R_{k,n} }
=
\esp{}{ \frac{L_{k,n}}{n} },
\end{equation}
but, for $1 \leq k \leq \varphi(n)$,
$
\esp{}{ L_{k,n}/{n} }
\leq
\esp{}{ B^n_{\varphi(n),1}/{n}}
$
and, as a consequence of
\cite[Theorem 1.1]{CL},
\[
\frac{B^n_{\varphi(n),1}}{n}
\build{\longrightarrow}{}{P}
0.
\]
Convergence of expectations follows, as
  $B^n_{\varphi(n),1}/n$ is bounded by 1.
\end{proof}

As a consequence, the  contribution of this regime is
\begin{equation}
\sum_{k=1}^{\varphi(n)}\esp{}{R_{k,n} }
\sim
\sum_{k=1}^{\varphi(n)}\frac{n}{n-k}
\sim
\sum_{k=n^{\frac{1}{2}+\varepsilon}}^{n-1}\frac{n}{k}
\sim
\pa{\frac{1}{2}-\varepsilon} n \log n.
\label{equationmpetit}
\end{equation}
\subsubsection{The transition regime.}
If  $k\simeq\sqrt n$,
$B^n_{k,\ell}=\Theta(n)$,
so that the  terms of the sum $R_n$ corresponding to
 the transition regime can be large.
However   there are few such terms:
\begin{eqnarray}
\label{equationmmoyen}
\sum_{k=\varphi(n)}^{\psi(n)}\esp{}{R_{k,n} }
\leq
\sum_{k=\varphi(n)}^{\psi(n)}
\frac{n}{n-k}
\sim
2 \varepsilon n \log n.
\end{eqnarray}

\subsubsection{The  almost full regime.}
If $k \geq \psi(n)$,
again as a consequence of
\cite[Theorem 1.1]{CL},
\begin{equation}
\label{bbb}
\frac{B^n_{\psi(n),1}}{n}
\build{\longrightarrow}{}{P}
1.
\end{equation}
Thus, as $L_{k,n}$ is the size of a size--biased pick among the blocks,
we expect that
$$
\pr{L_{k,n}\neq B^n_{k,1}}=o(1), \;
\frac{L_{k,n}}{n}\build{\longrightarrow}{}{P}1,
\mbox{ and }
\esp{}{R_{k,n} } = o\pa{\frac{n}{n-k}}.$$
More precisely, we have
\begin{lemm} \label{coucoumevoila_haut_regime}
$\displaystyle
\lim_n
\ \ \sup\ac{\left.
\frac{n-k}{n}\ \esp{}{R_{k,n} }\right|\ \psi(n) \leq k \leq n-1 }= 0.
$
\end{lemm}

\begin{proof}
Since  $L_{k,n}$ is the size of a size--biased pick among the blocks,
we should have
\[
\mathbb{P}(L_{k,n}=
B^n_{k,1}|B^n_{k,1})=\frac{B^n_{k,1}}{n},
\]
thus
$$
\esp{}{\frac{L_{k,n}}{n}}
\geq
\esp{}{\frac{L_{k,n}}{n} \mathbf{1}_{\{L_{k,n}=B^n_{k,1}\}}}
 \geq
\esp{}{\pa{\frac{B^n_{k,1}}{n}}^2}
 \geq
\esp{}{\pa{\frac{B^n_{\psi(n),1}}{n}}^2}.
$$
Now,
relations (\ref{aaa})
and
(\ref{bbb}) yields the desired result.
\end{proof}

Thus
\begin{equation}
\label{boundalmostfull}
\sum_{k=\psi(n)}^{n-1}\esp{}{R_{k,n}}
=
o\pa{\sum_{k=\psi(n)}^{n-1} \frac{ n}{n-k}}
=
o\pa{n \log n}.
\end{equation}

Lemma \ref{conv_en_moy_droit} follows, as (\ref{equationmpetit}), (\ref{equationmmoyen})
and (\ref{boundalmostfull}) hold true for any $\varepsilon$
positive and small enough.
\end{proof}

\begin{rema}
\label{Knuthmethod}
Note that, using
\begin{eqnarray*}
\esp{}{R_n}
& = &
\esp{}{ R_{n-1,n} } + \sum_{k=1}^{n-1} \frac{p_{n,k}+p_{n,n-k}}{2}
(\esp{}{R_{n-k} }+\esp{}{R_k }),
\end{eqnarray*}
and
\begin{eqnarray*}
\frac{p_{n,k}+p_{n,n-k}}{2}
& = &
\frac{1}{2(n-1)} C_n^k\left( \frac{k}{n}
\right)^{k-1}\left( \frac{n-k}{n} \right)^{n-k-1},
\end{eqnarray*}
we recover \cite[Relation (10.1)]{KS}. This lead Knuth and Sch{\"o}nhage
\cite{KS} to
 an alternative proof of Lemma
 \ref{conv_en_moy_droit}: one sees easily that
\[
\esp{}{ R_{n-1,n} }
=
a \sqrt n + O(1),
\]
in which $a=\sqrt{\pi/2}$, but
\cite[Relation (12.7)]{KS} ensures  that, as a consequence,
$$
\esp{}{C^{QFB}_{n,n-1}  }
=
\frac{a}{\sqrt{2 \pi}} n \log n + O(n).
$$
However, through this type of arguments,
we were  not able to obtain a suitable bound
for the variance.
\end{rema}

\subsection{Analysis of variance.}
\label{variance}

The next Proposition completes
the proof of Theorem \ref{olebotheo_right}.
\begin{prop}
\label{varianceR}
$
\var{C^{QFB}_{n,n-1}}
=
o((n\log n)^2).
$
\end{prop}
Once again, we use the exchangeability property of
 blocks' sizes in the parking scheme:
\begin{lemm}
\label{equation_ech_moy}
For $1 \leq l < k \leq n-1$,
$\displaystyle
\esp{}{R_{l,n}R_{k,n}}
=
\frac{\esp{}{R_{l,n}(n-L_{k,n})}}{n-k}.
$
\end{lemm}

\begin{proof}
Consider the $n-k+1$ blocks (clusters) before the $k$--th jump.
 Let us number them clockwise from 0 to $n-k$,
starting with the block  that contains the place filled by the
$l$--th arrival, and let $\gamma_i$ denote the size of the
$i$--th block. Let $Cl_0$ denote the random set of cars
belonging to block 0, and let $\mathcal F$ denote the $\sigma$--algebra
generated by $Cl_0$ and $\pa{t(c)}_{c\in Cl_0}$.
Also, let  $\mathcal G$  be the $\sigma$--algebra
generated by $\mathcal F$ and  $\pa{\gamma_i}_{1\le i\le n-k}$.
It is easy to see that,  among  the $(n-g_0)^{k-g_0-1}(n-k)$
possible parking configurations (given
$Cl_0$ and $\pa{t(c)}_{c\in Cl_0}$),  there are
\[
{k-g_0\choose g_1-1,g_2-1,\dots,g_{n-k}-1}\prod_{i=1}^{n-k}g_i^{g_i-2}
\]
configurations  such that
$\pa{\gamma_i}_{1\le i\le n-k}=\pa{g_i}_{1\le i\le n-k}$.
As a consequence, conditionally, given $\mathcal F$,
 the family $\pa{\gamma_i}_{1\le i\le n-k}$ is exchangeable,
while $\gamma_0$ is $\mathcal F$--measurable, and,
 being, in a sense, a size--biased pick among the $n-k+1$ blocks,
tends to be larger.
Note that $\gamma_0+...+\gamma_{n-k}=n$.

Given $\mathcal G$, the conditional probability that
 the $k$--th arrival fills the empty place at the end of block $i$
is $\gamma_{i}/n$, entailing that
\[
\esp{}{\left.R_{k,n}\right|\mathcal G}
=
 \frac{1}{n}
\sum_{i=0}^{n-k}\gamma_i
\gamma_{i+1},
\hspace{1 cm}
\esp{}{\left.L_{k,n}\right|\mathcal G}
=
 \frac{1}{n}
\sum_{i=0}^{n-k}\gamma_i^2,
\]
with the convention that ${n-k+1+\ell}={\ell}$.
As a consequence,
\[
\esp{}{R_{l,n}R_{k,n}}
=
\frac{1}{n}\ \esp{}{ R_{l,n}\sum_{i=0}^{n-k}\gamma_i\gamma_{i+1} }.
\]
Now, obviously, the relation
\[
\esp{}{ R_{l,n}\sum_{i=0}^{n-k}\gamma_i\gamma_{i+1}}
=
\esp{}{ R_{l,n} \sum_{i=0}^{n-k}
\gamma_{\sigma(i)}
\gamma_{\sigma(i+1)}  }
\]
holds when $\sigma $ is any power of the cyclic permutation
$(0,1,2,\dots,n-k)$, but, due to
the exchangeability of the sequence $(\gamma_i)_{1 \leq i \leq n-k+1}$,
 conditionally given $\mathcal F$, it also holds when $\sigma $
is any  permutation of the set $\{0,1,2,\dots,n-k\}$ leaving 0  invariant.
Thus, it holds for any $\sigma$, and, if $\mathfrak{S}_N$ is the set of  permutations on
$N$ elements:
\begin{eqnarray*}
\esp{}{R_{l,n}R_{k,n}}
& = &
\frac{1}{(n-k+1)!n} \sum_{\sigma \in \mathfrak{S}_{n-k+1}}
\esp{}{ R_{l,n} \sum_{i=0}^{n-k}
\gamma_{\sigma(i)}\gamma_{\sigma(i+1)}  }
\\
& = &
\frac{1}{(n-k)n}
\esp{}{ R_{l,n} \sum_{i=0}^{n-k}\sum_{j \neq i} \gamma_i\gamma_{j}}
\\
& = &
\frac{1}{n(n-k)}
\esp{}{ R_{l,n}\left( n^2 -\sum_{i=0}^{n-k} \gamma_i^2\right) },
\\
& = &
\frac{1}{n-k}
\esp{}{ R_{l,n}\esp{}{\left.n -L_{k,n}\right|\mathcal G} },
\end{eqnarray*}
completing the proof of the Lemma.
\end{proof}

Also, using the exchangeability property for the sequence
$(\beta_i)_{1 \leq i \leq n-k}$, as in Section \ref{variance}, we  obtain:
\begin{lemm} \label{equation_ech_carre}
For $1 \leq  k \leq n-1$,
$\displaystyle
\esp{}{R_{k,n}^2|L_{k,n}}
\leq
\frac{(n-L_{k,n})^2}{n-k}.
$
\end{lemm}

\begin{proof}[Proof of Proposition \ref{varianceR}.]
As in Section \ref{section:Meancost}, we decompose the variance
according to the three distinct regimes of the parking
scheme:
\begin{eqnarray*}
\lefteqn{\var{R_n}
 \leq
\sum_{k=1}^{n-1} \esp{}{R_{k,n}^2}
+
2\left( \sum_{1 \leq l <k \leq \varphi(n)}
\cov{R_{l,n},R_{k,n}} \right)
}
\\
&\quad &
+2\left( \sum_{{1 \leq l <k,}\atop{\varphi(n) \leq  k  \leq \psi(n)}}
  \esp{}{R_{l,n} R_{k,n}} \right)
+2\left( \sum_{{1 \leq l <k,}\atop{ \psi(n) \leq k \leq n}}
\esp{}{R_{l,n} R_{k,n}}\right).
\nonumber
\end{eqnarray*}

\subsubsection*{The square terms.}
By Lemma \ref{equation_ech_carre},  for  $1 \leq k \leq n-1$,
$\displaystyle
\esp{}{R_{k,n}^2}
\leq
\frac{n^2}{n-k},
$
so that
\begin{equation}
\label{equation_termescarres}
\sum_{k=1}^{n-1} \esp{}{R_{k,n}^2}
=
\bigO{n^2 \log n}.
\end{equation}

\subsubsection*{Covariances, the sparse regime.}
 Thanks to Lemma \ref{equation_ech_moy}, we have:
\begin{eqnarray}
\cov{R_{l,n}, R_{k,n}}
& \leq &
\esp{}{R_{l,n}}\left( \frac{n}{n-k} - \esp{}{R_{k,n}} \right).
\nonumber
\end{eqnarray}
This last inequality, combined with Lemma \ref{coucoumevoila},
entails that
\begin{eqnarray*}
\lim_n\ \ \sup_{ 1\le l<k \le \varphi(n)}\frac{(n-k)\cov{R_{l,n},  R_{k,n}}}{n\esp{}{R_{l,n}} }
&=&
0,
\end{eqnarray*}
so that
\begin{eqnarray}
\sum_{1 \leq l<k < \varphi(n)} \cov{R_{l,n}, R_{k,n}}
& = &
o\pa{
\sum_{l=1}^{\varphi(n)} \esp{}{R_{l,n}}
\sum_{k=l+1}^{\varphi(n)}\frac{n}{n-k}
}
\nonumber
\\
& = &
o\pa{
\sum_{l=1}^{\varphi(n)} \frac{n}{n-l}
\sum_{k=l+1}^{\varphi(n)}\frac{n}{n-k}
}
\nonumber
\\
&=&
o\pa{(n \log n)2}.
\label{equation_covpetits}
\end{eqnarray}

\subsubsection*{Covariances  when $k$ belongs to the
  transition regime. }
Thanks to Lemma \ref{equation_ech_moy}:
\begin{eqnarray}
\sum_{1 \leq l <k,\varphi(n) \leq  k  \leq \psi(n)}
\esp{}{R_{l,n} R_{k,n}}
& \leq & \sum_{1 \leq l \leq \psi(n)} \esp{}{R_{l,n}}
\sum_{\varphi(n) <  k\leq \psi(n)} \frac{n}{n-k}
\nonumber
\\
& \leq &
2 \varepsilon n \log n\sum_{1 \leq l \leq \psi(n)} \esp{}{R_{l,n}}
\nonumber
\\
& \leq &
2 \varepsilon (n \log n)^2.
\label{equation_croisesmoyens}
\end{eqnarray}

\subsubsection*{Covariances when $k$ belongs to the
  almost full regime. }
Note that  $\gamma=(\gamma_i)_{0\le i\le n-k}$ is the
 family of sizes of blocks before the $k$--th
arrival, numbered clockwise starting at some point that depends
on the $l$--th jump,
while $\beta=(\beta_i)_{0\le i\le n-k}$ is the same family,
numbered clockwise starting at some point that depends
on the $k$--th jump:
from the proof of Lemma \ref{equation_ech_moy},
we deduce that, for any $l<k$:
\begin{eqnarray*}
\sum_{l=1}^{k-1}\esp{}{R_{l,n} R_{k,n}}
&=&
\frac{1}{n(n-k)}
\esp{}{\sum_{l=1}^{k-1} R_{l,n}
\left( n^2 -\sum_{i=0}^{n-k} \beta_i^2\right) }.
\end{eqnarray*}
{From} expression (\ref{costadditivity}), we see that, conditionally,
 given that  $\beta=(b_i)_{0\le i\le n-k}$,
 the cost $\sum_{l=1}^{k-1} R_{l,n}$ is the sum of
$n-k+1$ random variables distributed as $(R_{b_i})_{0\le i\le n-k}$,
and, incidentally, independent. As a consequence of
Lemma \ref{conv_en_moy_droit}, there exists a universal constant $A$
such that
\begin{eqnarray*}
\esp{}{\sum_{l=1}^{k-1} R_{l,n}\vert \beta}
& \leq &
A
\esp{}{\sum_i \beta_i \log \beta_i}
\leq An\log n.
\end{eqnarray*}
Thus, for $k\ge\psi(n)$,
\begin{eqnarray*}
\sum_{l=1}^{k-1}\esp{}{R_{l,n} R_{k,n}}
& \leq &
\frac{A\log n}{n-k}\
\esp{}{ n^2 -\sum_{i=0}^{n-k} \beta_i^2}
\\
& \leq &
\frac{A\log n}{n-k}\
\esp{}{ n^2 -\max_{i} \beta_i^2}
\\
& \leq &
\frac{An^2\log n}{n-k}\
\esp{}{ 1 -\pa{\frac{B^n_{\psi(n),1}}{n}}^2}.
\end{eqnarray*}
Finally
\begin{eqnarray}
\label{equation_croisesgrands}
\sum_{1 \leq l <k,\psi(n) \leq  k  \leq n}
\esp{}{R_{l,n} R_{k,n}}
&\le&
o\pa{n^2\log n}
\sum_{\psi(n) \leq  k  \leq n} \frac{1}{n-k}
\end{eqnarray}

Again, since (\ref{equation_termescarres}),
(\ref{equation_covpetits}), (\ref{equation_croisesmoyens}) and
(\ref{equation_croisesgrands}) hold true for any $\varepsilon$
positive and small enough,
this completes the proof of Proposition \ref{varianceR}.
\end{proof}

\begin{rema}
While the  asymptotic behaviour of the partial costs was
 obtained by merely analytic tools, our analysis
of the complete costs
relies on the additional
 information captured by some underlying
combinatorial structure,  the parking scheme,
and can hardly be extended to other kernels.
\end{rema}

\section{Asymptotics of the cost of Quick Find}
\label{QFfullsection}

This Section is devoted to the proof of  Theorem
\ref{olebotheo_random}. We need some notations.
First, as the cost $A_{k,n}$  of the $k$--th
union of a Quick Find algorithm
is a random uniform pick among the sizes of the
two clusters involved, we may write
$$A_{k,n}=\varepsilon_k L_{k,n}+(1-\varepsilon_k) R_{k,n},$$
in which
$(\varepsilon_k)_{1 \leq k \leq n-1}$
is a sequence of i.i.d. random variables
with law
$\frac{1}{2}\delta_0+\frac{1}{2}\delta_1$, independent of the parking
scheme.
Also, let $c(k)$ denote the car involved in the $k$--th
jump, that is, such that
\[\#\ac{c\left\vert1\le c\le n-1\mbox{ and }T_c\le T_{c(k)}\right.}
=
k,
\]
let  the first try of $c(k)$, $t(c(k))$, be denoted $t(k)$ for sake of  brevity,
and let $f(k)$ be the final place of
of $c(k)$.
Let $\mathcal H$  (resp. $\mathcal H_k$)
be the $\sigma$--algebra
generated by  $\pa{t(c),T_c}_{c\in \mathcal C}$
(resp. by $\pa{t(i)}_{1\le i\le k-1}$  and $f(k)$).
Finally, set
\begin{eqnarray*}
 F_{k,n}
&=&
\frac{1}{2}(L_{k,n}+R_{k,n}),
\\
F_n
=
\sum_{k=1}^{n-1} F_{k,n},
\mbox{ }\mbox{ }
L_n
&=&
\sum_{k=1}^{n-1} L_{k,n},
\mbox{ }\mbox{ }
 D_n
=
\sum_{k=1}^{n-1} D_{k,n}
.
\end{eqnarray*}
The proof is based on  the following observations: clearly
\begin{equation}
\label{rem_cout_cond}
\esp{}{A_{k,n}|\mathcal H}
=
\frac{1}{2}(L_{k,n}+R_{k,n}),
\end{equation}
and, since, conditionally given $L_{k,n}$, the displacement $D_{k,n}$ is
uniformly distributed
on $\{1, ..., L_{k,n}\}$, we have
\begin{equation}
\label{rem_depl_cond}
\esp{}{D_{k,n}|
\mathcal H_k}
=
\frac{1}{2}(L_{k,n}+1).
\end{equation}
We also need an important result about hashing with linear probing
\cite{CM,FPV,Ja}:
\begin{theo}[Flajolet, Poblete and Viola, 1998]
$$n^{-3/2}{D_n} \build{\longrightarrow}{}{law} \int_0^1
e(t)dt.$$
\end{theo}

Due to relation (\ref{rem_depl_cond}), we have
\begin{lemm}
\label{lemme1}
$\displaystyle
\left\Vert2D_n-L_{n}\right\Vert_2
=
o\pa{n^{3/2}}.
$
\end{lemm}

\begin{proof}
Expanding
$(2D_n-L_{n}-n+1)^2$,
we obtain:
\[
\left\Vert2D_n-L_{n}-n+1\right\Vert_2^2
=
\Xi_1 +\Xi_2,
\]
in which
\begin{eqnarray}
\Xi_1
& = &
\sum_{k=1}^n\ \esp{}{\pa{2D_{k,n}-L_{k,n}-1}^2 }
\nonumber
\\
\Xi_2
& = &
2 \sum_{1\le i<j\le n-1} \esp{}{ \left(2D_{i,n}-L_{i,n}-1
\right)\left(2D_{j,n}-L_{j,n}-1\right)}.
\nonumber
\end{eqnarray}
Owing to (\ref{rem_depl_cond}), for  $i<j$,
\begin{eqnarray*}
\esp{}{
\esp{}{ \left(2D_{i,n}-L_{i,n}-1\right)\left(2D_{j,n}-L_{j,n}-1\right)
\left| \
\mathcal H_j
\right.}}
& = & 0,
\end{eqnarray*}
and $\Xi_2$ vanishes.
By definition of $D_{k,n}$, we also have
\begin{eqnarray}
\esp{}{ \left.\left( 2D_{k,n}-L_{k,n}-1
 \right)^2 \right| L_{k,n}}
& = & \frac{1}{3}\pa{L_{k,n}^2-1}
\nonumber
\end{eqnarray}
Thus
\begin{eqnarray}
\Xi_1
&\leq&
  \frac{1}{3}
\sum_{k=1}^{n-1} \esp{}{L_{k,n}^2}
\le
\frac{n^3}{3} \int_0^1 \esp{}{\left( \frac{L_{\lceil\alpha  n\rceil,n}}{n}
    \right)^2 } d \alpha.
 \label{equ_conv_int}
\end{eqnarray}
\noindent
According to \cite{PITT}, for $0< \alpha <1$,
 $( B^n_{\lceil\alpha n\rceil,1}/n)_{n \in \mathbb{N}}$ converges in
probability to $0$, thus
\[
\lim_n\esp{}{ \pa{\frac{L_{\lceil\alpha n\rceil,n}}n}^2}
=
0
\]
and Lebesgue Dominated
Convergence Theorem completes the proof.
\end{proof}

As a consequence of Lemma \ref{conv_en_moy_droit}
and  Proposition \ref{varianceR},
\begin{lemm}
\label{lemme2}
$\displaystyle
\left\Vert2F_n-L_{n}\right\Vert_2
=
\left\Vert R_n\right\Vert_2
=
o\pa{n^{3/2}}.
$
\end{lemm}

Finally,
\begin{lemm} \label{lemme3}
$\displaystyle
\left\Vert F_n-C^{QF}_{n,n-1}\right\Vert_2
=
o\pa{n^{3/2}}.
$
\end{lemm}

\begin{proof}
We split
\begin{eqnarray*}
\left\Vert F_n-C^{QF}_{n,n-1}\right\Vert_2^2
 & = &
\esp{}{ \left( \sum_{k=1}^{n-1}
\left(\varepsilon_k -\frac{1}{2} \right) L_{k,n}
+\left(\frac{1}{2}-\varepsilon_k \right) R_{k,n} \right)^2 }
\nonumber
\end{eqnarray*}
in three terms:
\begin{eqnarray*}
\Xi_1
& = &
\esp{}{\left( \sum_{k=1}^{n-1}\left(\varepsilon_k -\frac{1}{2} \right)
L_{k,n} \right)^2 }
\\
\Xi_2
& = &
\esp{}{ \left( \sum_{k=1}^{n-1}
\left( \frac{1}{2}-\varepsilon_k \right) R_{k,n} \right)^2 }
\\
\Xi_3
& = &
2\ \sum_{i,j} \esp{}{\left(\varepsilon_i -\frac{1}{2} \right)
\left(\frac{1}{2}-\varepsilon_j
  \right)L_{i,n} R_{j,n} }
\end{eqnarray*}
Since
$\left(\varepsilon_k -\frac{1}{2}\right)_{1 \leq k \leq n-1}$
are i.i.d. random variables with mean $0$,
independent of $\mathcal H$,  we find,
conditioning to $\mathcal H$,  that:
\begin{eqnarray*}
\Xi_1
&=&
\frac{1}{4} \sum_{k=1}^{n-1} \esp{}{L_{k,n}^2},
\\
\Xi_2
&=&
\frac{1}{4} \sum_{k=1}^{n-1} \esp{}{R_{k,n}^2},
\\
\Xi_3
&=&
-\frac{1}{2} \sum_{k=1}^{n-1} \esp{}{ L_{k,n} R_{k,n}}.
\end{eqnarray*}
We conclude using the same arguments as  in the proof of   (\ref{equ_conv_int}),
since we have
\begin{eqnarray*}
\left\Vert F_n-C^{QF}_{n,n-1}\right\Vert_2^2
& = &
\frac{1}{4} \sum_{k=1}^{n-1} \esp{}{\left(L_{k,n} -R_{k,n}
  \right)^2 }
\\
& \leq &
\frac{n3}{4} \int_0^1
\esp{}{\left( \frac{L_{\lceil\alpha n\rceil,n}+R_{\lceil\alpha  n\rceil,n}}{n}
    \right)^2 } d \alpha.
\end{eqnarray*}
\end{proof}

Finally Theorem \ref{olebotheo_random}
is obtained by
combining these Lemmas with
 \cite[Theorem 4.1]{Bi}:
\begin{theo}
\label{thetrick}
Let $(X_n)_{n \in \mathbb{N}}$, $(Y_n)_{n \in \mathbb{N}}$
 and $X$ be random variables
such that for every $n$, $X_n$ and $Y_n$ are defined on the same
probability space. If $(X_n)_{n \in \mathbb{N}}$ converges in law to  $X$ and if
$(|X_n-Y_n|)_{n \in \mathbb{N}}$ converge in probability to $0$ then  $(Y_n)_{n \in
  \mathbb{N}}$ converges in law to $X$.
\end{theo}

\section{Almost full regime: Proof of Theorem \ref{process}}
\label{almostfullsection}

Here we list the slight adaptations
  to be made to the previous proof, in order to obtain Theorem
\ref{process}.
We introduce
\[
D_n(\beta)
=
n^{-3/2}\sum_{k=1}^{\lfloor n-\beta\sqrt n\rfloor}
D_{k,n}
\]
and we observe that by the same proof
as in the previous Section, but considering
partial sums rather than
the complete sums, we obtain
\begin{equation}
\label{dispapprox}
\left\Vert D_n(\beta)
-W_n(\beta)\right\Vert_2^2
=
o(1).
\end{equation}
On the other hand, as a direct consequence of \cite {CL} (see specially
  \cite[Theorem 4.1]{CL}), we know that   it is possible to build,
on a suitably chosen
  probability space $\Omega$, a version of the normalized
Brownian excursion, and also
a version of the parking scheme for each possible size $m$,
in such a way that,
  if  $\sqrt m\ \psi_m\left(\beta,t\right)$ denotes
  the number of cars that tried to park, successfully or not,
  on place $\lfloor tm\rfloor$, among the $\lfloor m-\beta\sqrt  m\rfloor$
  cars already arrived,
then we have:
\[
\Pr\left(\forall \Lambda,
\hspace{0.2cm}
\psi_m(\beta,t)
\hspace{0,3cm}
\build{\tend}{on\textrm{ }\Delta_{\Lambda}}{uniformly}
\hspace{0,3cm}
h_{\lambda}(t)\right)=1,
\]
in which $\Delta_{\Lambda}=[0,\Lambda]\times[0,1]$.

Since
$\psi_m$ captures the whole story of the parking process
(for instance, it captures
the sizes and positions of blocks and the
first tries of successive cars),
$\psi_m$ also describes the sample
paths of the additive Marcus--Lushnikov
processes with size $m$. Specifically, the total and partial
displacements
have the following simple expression in terms of $\psi_m$:
\[
D_n(\beta)
=
\int_0^1\psi_m\pa{\beta,t}\,dt.
\]
{From} this relation, we obtain directly that
\[
\Pr\left(\forall \Lambda,
\hspace{0.2cm}
D_n(\beta)
\hspace{0,3cm}
\build{\tend}{on\textrm{ }[0,\Lambda]}{uniformly}
\hspace{0,3cm}
W(\beta)\right)=1,
\]
which, together with
\eqref{dispapprox},
entails the convergence of finite--dimensional distributions of
the positive decreasing  processes
$W_n(\cdot)$  to the finite--dimensional distributions of
$W(\cdot)$. This is enough to insure the weak convergence of these
processes,
seen as random variables with values in the space of
tail distributions of positive measures on $[0,+\infty]$, endowed
with the topology of weak convergence of the corresponding
positive measures.
  These spaces are Lusin spaces, thus,
according to  the Skorohod representation theorem  \cite[II.86.1]{ROW},
one can find a probability space where the weak convergence of
$W_n(\cdot)$  to
$W(\cdot)$ is almost sure
and since $\beta\rightarrow W(\beta)$ is almost surely
continuous, it entails that
$W_n(\cdot)$  converges to
$W(\cdot)$ uniformly on $[0,+\infty]$,
almost surely on the
  probability space $\Omega$.

\section{Concluding remarks}

Knuth and Sch{\"o}nhage gave asymptotics
for the \textit{expectation} of some additive functionals
of the additive Marcus--Lushnikov process, and
  we were able to give a more precise information,   either
  the asymptotic behaviour
of the \textit{distribution}, or a concentration result, for
these functionals, by embedding the additive
Marcus--Lushnikov process in a richer structure.
It would be interesting to extend such results to
     Marcus--Lushnikov processes
with a general kernel $K(x,y)$, but general theorems of convergence
of Marcus--Lushnikov processes seem not precise enough,
at least for the total costs,
to allow such a generalisation right now. For the total costs, our
  approach is quite specific of the additive case,
and even in the important case
$K(x,y)=xy$ it seems rather hard to improve the results of
  Bollob{\'a}s \&  Simon \cite{BOLL},
  who show that
the average cost of QFW  is $cn+O(n/\log n)$, $c=2.0847\cdots$,
  while the average of QF is $n^2/8+O\pa{n(\log n)^2}$.

\section*{Acknowledgements}
We would like to thank Philippe Flajolet for pointing to us this  problem,
in relation with the ``Cutting down random trees"
problem of Meir \& Moon \cite{mmoon},
that we learned from Jean-Fran\c cois  Marckert.
The second author also thanks Nicolas Fournier for many fruitful
discussions. 

\def\refname{References}

\end{document}